\documentclass[a4paper,12pt,reqno]{amsart}
\usepackage[T1]{fontenc}
\usepackage[utf8]{inputenc}
\usepackage{amsthm,amsmath,amssymb,mathrsfs}
\usepackage{graphicx}
\usepackage{xcolor}
\usepackage[colorlinks=true, pdfstartview=FitV, linkcolor=blue, citecolor=red, urlcolor=blue]{hyperref}

\numberwithin{equation}{section}
\newtheorem{theorem}{Theorem}[section]

\newtheorem{definition}[theorem]{Definition}
\newtheorem{example}[theorem]{Example}
\newtheorem{lemma}[theorem]{Lemma}
\newtheorem{proposition}[theorem]{Proposition}
\newtheorem{remark}[theorem]{Remark}

\allowdisplaybreaks

\newcommand{\RR}{\mathbb{R}}

\DeclareMathOperator*{\esup}{ess\,sup}
\DeclareMathOperator*{\einf}{ess\,inf}

\title[Weighted $H^p$--$L^q$ estimates for operators with multiple kernels]{Weighted $H^p$--$L^q$ estimates for fractional operators with multiple kernels under a strengthened H\"ormander condition}

\author[A. L. Gallo]{Andrea L. Gallo}
\address{A. L. Gallo\\ FaMAF \\ Universidad Nacional de C\'ordoba \\ CIEM (CONICET) \\ 5000 C\'ordoba, Argentina}
\email{andrea.gallo@unc.edu.ar}

\author[M. S. Riveros]{M. Silvina Riveros}
\address{M. S. Riveros\\ FaMAF \\ Universidad Nacional de C\'ordoba \\ CIEM (CONICET) \\ 5000 C\'ordoba, Argentina}
\email{msriveros@unc.edu.ar}

\author[L. A. Vallejos]{Lucas A. Vallejos}
\address{L. A. Vallejos\\ FaMAF \\ Universidad Nacional de C\'ordoba \\ CIEM (CONICET) \\ 5000 C\'ordoba, Argentina}
\email{lucas.vallejos@unc.edu.ar}

\thanks{The authors are partially supported by CONICET and SECYT-UNC}
\date{}

\begin{document}
	
	\begin{abstract}
		We study a scale of strengthened fractional H\"ormander conditions, indexed
		by $\delta\ge0$ and adapted to the Hardy spaces $H^p$, $0<p<1$, which was
		recently considered by Han and Wu [Y.~Han, H.~Wu, \emph{The weighted $H^p$
			estimates for commutators of fractional integrals}, Potential Anal.
		\textbf{59} (2023), 1827--1850] for fractional integrals with a single
		kernel. We show that these conditions form a strictly decreasing scale in
		$\delta$, and we compute the exact range of $\delta$ for a natural family of
		homogeneous kernels. We then apply them to the integral operators
		\[
		T_\alpha f(x)=\int_{\RR^{n}}k_1(x-A_1y)\cdots k_m(x-A_my)\, f(y)\,dy,
		\]
		where $A_1,\dots,A_m$ are invertible matrices and
		$\frac n{q_1}+\dots+\frac n{q_m}=n-\alpha$, $0<\alpha<n$. Assuming that each
		kernel $k_i$ satisfies an $L^{p_i}$ size condition of order $n-n/q_i$ and the
		H\"ormander condition of order $\delta=n(\frac1p-1)$, we obtain, for
		$\max\{\frac n{n+1},\frac n{n+\alpha}\}<p<1$ and $1/q=1/p-\alpha/n$, the
		$H^{p}_{w^p}(\RR^{n})-L^{q}_{w^q}(\RR^{n})$ boundedness of $T_\alpha$ for
		weights $w$ such that $w^{a}\in\mathcal A_1$, where
		$\frac1a=\frac1q-\sum_i\frac1{p_i}$, and $w(A_ix)\le c\,w(x)$ a.e. For a
		single kernel, this relaxes the smoothness and weight conditions of a recent
		result of Han and Wu. No homogeneity is assumed on the kernels; the case
		$k_i(x)=\Omega_i(x)|x|^{-n/q_i}$ with $\Omega_i$ homogeneous of degree zero is
		included. Finally, we show that the strengthened condition cannot be
		replaced by the classical fractional H\"ormander condition.
	\end{abstract}
	\maketitle
	
	\section{Introduction}
	
	The Hardy spaces on $\RR^n$ were introduced in \cite{FS} by C. Fefferman
	and E. Stein, and have been extensively studied ever since. A key feature of
	these spaces is that they provide an effective alternative to Lebesgue spaces
	in the range $0<p<1$. In particular, while Riesz transforms fail to be
	bounded on $L^p(\RR^n)$ for $p\leq 1$, they are bounded on the
	corresponding Hardy spaces $H^p(\RR^n)$. To investigate the boundedness of
	operators such as singular integrals or fractional-type operators on Hardy
	spaces, one typically relies on their atomic or molecular characterizations.
	The atomic decomposition was first established by Coifman in one dimension
	(see \cite{C}) and later extended to higher dimensions by Latter \cite{L}.
	
	These descriptions suggest that the boundedness of linear operators on $H^p$
	could, in principle, be reduced to analyzing their action on atoms. However,
	this approach has certain limitations. Bownik \cite{B}, building on an
	example due to Meyer, constructed a linear functional defined on a dense
	subspace of $H^1(\RR^n)$ that maps all $(1,\infty,0)$-atoms into bounded
	scalars, but does not admit a bounded extension to $H^1(\RR^n)$.
	Nevertheless, for classical operators $T$, uniform boundedness on atoms
	implies boundedness from $H^p$ into $L^p$. This follows from the fact that
	$T$ is bounded on $L^s$ for some $1<s<\infty$, together with the possibility
	of choosing atomic decompositions that converge in $L^s$; see
	\cite{YZ} and \cite{Rocha-2}.
	
	The weighted Hardy spaces $H^p_{w}(\RR^n)$ are defined by adapting the
	definition of $H^p(\RR^n)$ to the measure $w(x)\,dx$ (see \cite{ST}), and
	the theory is particularly meaningful for $w\in\mathcal A_\infty$. The
	atomic decomposition of $H^p_{w}(\RR^n)$ was established in \cite{GC,ST}.
	In \cite{Rocha}, the author obtained, for all $w\in\mathcal A_\infty$, an
	atomic decomposition of $H^p_w(\RR^n)$ converging in $L^s$, $1<s<\infty$,
	and a molecular decomposition. As applications, he obtained the boundedness
	of certain singular integrals on $H^p_w$ and from $H^p_w$ into $L^p_w$, for
	all $w\in\mathcal A_\infty$ and $0<p\le1$, and the $H^p_{w^p}-H^q_{w^q}$ and
	$H^p_{w^p}-L^q_{w^q}$ boundedness of the Riesz potential $I_\alpha$, for
	$\frac1q=\frac1p-\frac\alpha n$ and certain weights $w$.
	
	Let $\Omega\in L^r(S^{n-1})$, $r>1$, be homogeneous of degree zero, and
	consider
	\[
	T_{\Omega,\alpha}f(x)=\int\frac{\Omega(x-y)}{|x-y|^{n-\alpha}}f(y)\,dy.
	\]
	Weighted estimates for $T_{\Omega,0}$ were obtained in \cite{KuW,Du,Wa},
	and for $\alpha>0$ in \cite{DL}. In \cite{BLRi}, the authors obtained a
	Coifman type inequality for general fractional integral operators whose
	kernels satisfy a size condition and a H\"ormander condition given by a
	Young function. In \cite{DL2}, Ding and Lu proved, under $L^r$-Dini type
	conditions on $\Omega$, the $H^p-L^q$ boundedness of $T_{\Omega,\alpha}$ for
	$p=1$ and, when $0<\alpha<1$, for $\frac n{n+\alpha}\le p<1$, as well as its
	$H^p-H^q$ boundedness for a range of $p<1$ when $0<\alpha<1/2$.
	Str\"omberg and Wheeden \cite{SW} gave the weighted $H^p-L^q$ and
	$H^p-H^q$ boundedness of $I_\alpha$, and Ding, Lee and Lin \cite{DLeLi},
	using the atom--molecule theory of \cite{GC,LeLi}, obtained the weighted
	$H^p-L^q$ and $H^p-H^q$ boundedness of $T_{\Omega,\alpha}$. More recently,
	Han and Wu \cite{HW} considered fractional H\"ormander conditions of order
	$\delta$, which strengthen those of \cite{BLRi}, and obtained weighted
	$H^p-L^q$ estimates for fractional integrals $\int K_\alpha(x-y)f(y)\,dy$
	with a single kernel satisfying such a condition; see
	\cite[Definition~5.6 and Theorem~5.10]{HW}. Estimates for $I_\alpha$ and for
	$T_{\Omega,\alpha}$ with $\Omega$ Lipschitz on Musielak--Orlicz Hardy spaces
	were recently obtained in \cite{HHSW}.
	
	Let $0\le\alpha<n$ and $m\in\mathbb N$. For $1\le i\le m$ let
	$1<q_i<\infty$ be such that
	\begin{equation}\label{q_i}
		\frac n{q_1}+\dots+\frac n{q_m}=n-\alpha .
	\end{equation}
	We study the integral operator
	\begin{equation}\label{operador}
		T_\alpha f(x)=\int_{\RR^{n}}K(x,y)f(y)\,dy,\qquad
		K(x,y)=k_1(x-A_1y)\cdots k_m(x-A_my),
	\end{equation}
	where $A_1,\dots,A_m$ are invertible matrices, $k_1,\dots,k_m$ are
	measurable functions on $\RR^n$, and $f\in L^\infty_c(\RR^n)$.
	
	When $A_i=a_iI$ and $k_i(x)=|x|^{-n/q_i}$, Godoy and Urciuolo \cite{GU2}
	obtained the $L^p-L^q$ boundedness of $T_\alpha$ for $1<p<\frac n\alpha$,
	$\frac1q=\frac1p-\frac\alpha n$. For smooth $\Omega_i$, Rocha and Urciuolo
	\cite{RoU} considered matrices satisfying
	
	\smallskip
	(H)\ \ \emph{$A_{i}$ is invertible and $A_i-A_j$ is invertible for
		$i\neq j$, $1\leq i,j\leq m$,}
	\smallskip
	
	\noindent and proved that $T_\alpha$ is bounded from $H^{p}(\RR^{n})$ into
	$L^{q}(\RR^{n})$ for $0<p<\frac n\alpha$; they also showed that
	$H^p-H^q$ boundedness cannot be expected, an important difference with the
	case $m=1$. The weighted $L^p-L^q$ boundedness of $T_\alpha$ with rough
	homogeneous kernels was obtained by Riveros and Urciuolo \cite{RiU2}, and
	Ib\'a\~nez Firnkorn and Riveros \cite{FR} treated kernels satisfying size
	and H\"ormander conditions of Young type. More recently, Wu and Zhang
	\cite{WZ} studied these operators, with rough homogeneous kernels, on
	Orlicz--Morrey spaces.
	
	In this paper we do not assume that the kernels are homogeneous. Following
	\cite[Definition~3.1]{BLRi}, we only require a size condition on dyadic
	annuli. For the regularity of the kernels we use the fractional
	H\"ormander conditions of order $\delta$ of \cite{HW}, which strengthen the
	fractional H\"ormander condition of \cite[Definition~3.2]{BLRi} and are
	adapted to the Hardy spaces $H^p$ with $p<1$. We study the structure of this
	scale of conditions and its sharpness, and we use it to obtain weighted
	$H^p-L^q$ estimates for operators with several kernels composed with linear
	maps. These are the main contributions of this paper.
	
	\subsection*{Hypotheses and main result}
	We write $|x|\sim s$ for the annulus $\{x\in\RR^n: s<|x|\le 2s\}$ and, for
	$1\le\tau<\infty$ and $s>0$, we consider the normalized annular norm
	\[
	\|f\|_{\tau,|x|\sim s}
	=\Bigl(\frac1{|B(0,2s)|}\int_{s<|x|\le2s}|f(x)|^\tau\,dx\Bigr)^{1/\tau},
	\]
	with $\|f\|_{\infty,|x|\sim s}=\esup_{s<|x|\le 2s}|f(x)|$; this is the norm
	$\|f\|_{A,|x|\sim s}$ of \cite{BLRi} with $A(t)=t^\tau$. Throughout the
	paper $0<\alpha<n$, $0<p<1$,
	\[
	\delta=n\Bigl(\frac1p-1\Bigr),\qquad \frac1q=\frac1p-\frac\alpha n .
	\]
	Besides (H) (which, for $m=1$, only requires $A_1$ to be invertible), we
	assume:
	
	\smallskip
	($\mathrm{H}_1$)\ \ For each $1\le i\le m$ there is $q_i<p_i\le\infty$
	such that $\|k_i\|_{p_i,|x|\sim s}\le C s^{-n/q_i}$ for all $s>0$.
	
	\smallskip
	($\mathrm{H}_2$)\ \ For each $1\le i\le m$, with the same $p_i$ as in
	$(\mathrm H_1)$, there exist $c_i\ge1$ and $C_i>0$ such that
	\[
	\sum_{j=1}^{\infty}2^{j\delta}(2^jR)^{n/q_i}
	\|k_i(\cdot-h)-k_i\|_{p_i,|x|\sim2^jR}\le C_i,
	\qquad h\in\RR^n,\ R>c_i|h|.
	\]
	
	\smallskip
	\noindent Condition $(\mathrm H_1)$ is the size condition of
	\cite[Definition~3.1]{BLRi}. Condition $(\mathrm H_2)$ means that
	$k_i\in H^{(\delta)}_{n-n/q_i,p_i}$, where, for $1\le\tau\le\infty$,
	$0\le\beta<n$ and $\delta\ge0$, $H^{(\delta)}_{\beta,\tau}$ denotes the class
	of kernels satisfying the fractional H\"ormander condition of order $\delta$
	introduced by Han and Wu \cite[Definition~5.6]{HW} (see
	Definition~\ref{def-Hdelta}). For $\delta=0$ we recover the
	fractional H\"ormander condition of \cite[Definition~3.2]{BLRi} (the case
	$k=0$ there), while for $\delta=n(\frac1p-1)$ the factor $2^{j\delta}$ is
	exactly the one produced by the action of $T_\alpha$ on the atoms of
	$H^p_{w^p}$ (see Remark~\ref{R-origin}).
	
	The classes $H^{(\delta)}_{\beta,\tau}$ form a genuine scale. They decrease
	as $\delta$ increases, and Proposition~\ref{P-strict} shows that, for
	$n\ge2$, all these inclusions are strict:
	\[
	S_{\beta,\tau}\cap H^{(\delta_2)}_{\beta,\tau}\subsetneq
	S_{\beta,\tau}\cap H^{(\delta_1)}_{\beta,\tau}
	\qquad\text{for } 0\le\delta_1<\delta_2,\ \delta_1<1,\ 1\le\tau\le\infty .
	\]
	Moreover, the position of a kernel in this scale can be computed explicitly
	in natural examples: for the homogeneous kernels
	$k_\varepsilon(x)=|x_n|^{\varepsilon}|x|^{\beta-n-\varepsilon}$, with
	$-\frac1\tau<\varepsilon<1-\frac1\tau$ and $\varepsilon\neq0$, we have
	$k_\varepsilon\in H^{(\delta)}_{\beta,\tau}$ if and only if
	$\delta<\varepsilon+\frac1\tau$. In particular, for every $\delta>0$
	condition $(\mathrm H_2)$ is strictly stronger than the classical one, and it
	becomes more restrictive as $p$ decreases. These classes are studied in
	Section~\ref{S-Hormander}.
	
	Since $p_i>q_i$, \eqref{q_i} gives
	$\sum_{i}\frac1{p_i}<\sum_i\frac1{q_i}=1-\frac\alpha n<1$. We define
	$t\ge1$ and $a=a(t,q)$ by
	\begin{equation}\label{def-a}
		\frac1{p_1}+\dots+\frac1{p_m}+\frac1t=1,
		\qquad
		\frac1{a}=\frac1q-\sum_{i=1}^m\frac1{p_i}=\frac1q+\frac1t-1,
	\end{equation}
	that is, $a(t,q)=\frac{tq}{t+q-tq}$. Since $p<1$, we have
	$\frac1q=\frac1p-\frac\alpha n>1-\frac\alpha n>\sum_i\frac1{p_i}$, so $a$ is
	well defined and positive. Note also that $t<\frac n\alpha$ and $a\ge q$,
	with equality only if $p_i=\infty$ for every $i$.
	
	The main result of this paper is the following.
	
	\begin{theorem}\label{T-HpLq}
		Let $0<\alpha<n$ and $\max\Bigl\{\frac n{n+1},\frac n{n+\alpha}\Bigr\}<p<1, \ \delta=n\Bigl(\frac1p-1\Bigr), \ \frac1q=\frac1p-\frac\alpha n .$
		
		Let $T_\alpha$ be defined by \eqref{operador}, and assume $(\mathrm H)$,
		$(\mathrm H_1)$ and $(\mathrm H_2)$. If $w$ is a weight such that
		$w^{a(t,q)}\in\mathcal A_1$ and $w(A_ix)\le C\,w(x)\ \text{ a.e. } x\in\RR^n,\ 1\le i\le m,$
		then $T_\alpha$ extends to a bounded operator from $H^p_{w^p}(\RR^n)$ into
		$L^q_{w^q}(\RR^n)$.
	\end{theorem}
	
	The restriction on $p$ is equivalent to $\delta<1$ and $q>1$. The
	condition $\delta<1$ allows us to use atoms with vanishing integral only
	(see the proof of Theorem~\ref{T-HpLq}), and $q>1$ is used in Minkowski's
	integral inequality in Lemma~\ref{TL}. Moreover, $\delta<1$ is natural for
	the hypotheses themselves: if $\delta>1$, then $(\mathrm H_1)$ and
	$(\mathrm H_2)$ force $k_i\equiv0$ (Remark~\ref{R-delta-large}).
	
	Let us comment on the hypotheses of Theorem~\ref{T-HpLq}. Although
	$(\mathrm H_2)$ is strictly stronger than the classical fractional
	H\"ormander condition, the classes $H^{(\delta)}_{\beta,\tau}$ are large:
	conditions $(\mathrm H_1)$ and $(\mathrm H_2)$ are satisfied by many
	kernels, homogeneous or not, such as kernels with pointwise H\"older regularity of
	order $\varepsilon>\delta$, rough homogeneous kernels
	$\Omega_i(x)|x|^{-n/q_i}$ satisfying an $L^{p_i}$-Dini condition of order
	$\delta$, and kernels obtained from Calder\'on--Zygmund kernels (see
	Section~\ref{S-Hormander}). For homogeneous kernels, the extra decay required
	in $(\mathrm H_2)$ cannot be obtained by changing the index of a classical
	fractional H\"ormander class (Theorem~\ref{T-index}), and the factor
	$2^{j\delta}$ cannot be removed: in Section~\ref{S-sharp} we give an
	operator whose kernels satisfy $(\mathrm H_1)$ and the classical fractional
	H\"ormander condition, but for which the conclusion of
	Theorem~\ref{T-HpLq} fails; in fact, in general $\delta$ cannot be replaced
	in $(\mathrm H_2)$ by any $\delta'<\delta-\frac1q+\frac1{p_1}$
	(Remark~\ref{R-counterexample}). Theorem~\ref{T-HpLq} also covers the case
	$m=1$ of a single kernel. In that case, Han and Wu \cite[Theorem~5.10]{HW}
	obtained a related weighted estimate, for $0<\alpha\le1$, assuming the
	H\"ormander condition of order $\alpha$; Theorem~\ref{T-HpLq} only requires
	the condition of order $\delta=n(\frac1p-1)\le\alpha$, allows any
	$0<\alpha<n$, and has a weaker condition on the weight (see
	Remark~\ref{R-HW} for a detailed comparison). Finally,
	the hypotheses on $w$ and on the matrices are satisfied, for instance, by the power weights
	$w(x)=|x|^{\gamma}$, $-n/a(t,q)<\gamma\le0$, and arbitrary matrices satisfying
	$(\mathrm H)$; more examples are given in Section~\ref{S-examples}.
	
	The paper is organized as follows. In Section~\ref{S-prelim} we recall some
	preliminaries on weights and weighted Hardy spaces.
	Section~\ref{S-Hormander} is devoted to the classes
	$H^{(\delta)}_{\beta,\tau}$: we recall their definition, prove that they form a
	strictly decreasing scale (Proposition~\ref{P-strict}), show that
	$0\le\delta<1$ is the only meaningful range (Remark~\ref{R-delta-large}),
	give sufficient conditions covering kernels with pointwise H\"older
	regularity, rough homogeneous kernels with Dini conditions of order
	$\delta$ and kernels obtained from Calder\'on--Zygmund kernels, and study
	the homogeneous case (Theorem~\ref{T-index}). Section~\ref{S-proof} contains
	the proof of Theorem~\ref{T-HpLq}, which is built on these classes: an
	extension theorem with two weights, the technical lemma, where
	$(\mathrm H_2)$ is used precisely to absorb the factor $2^{j\delta}$
	produced by the atoms, the unweighted estimate, the proof itself, a
	comparison with \cite{HW}, and the sharpness of $(\mathrm H_2)$. In Section~\ref{S-examples} we
	give examples of matrices and weights satisfying the hypotheses of
	Theorem~\ref{T-HpLq}.
	
	\section{Preliminaries}\label{S-prelim}
	
	\subsection{Weights}
	A weight is a function $w\in L^1_{loc}(\RR^n)$ with values in $(0,\infty)$
	almost everywhere. Given a weight $w$ and $0<p<\infty$, $L^{p}_w(\RR^{n})$
	is the space of measurable functions with
	$\|f\|^{p}_{L^{p}_w}=\int_{\RR^{n}}|f|^{p}w<\infty$. For a measurable set
	$E$ and $0<p<\infty$ we write $w^p(E)=\int_E w^p$. For $p=\infty$,
	$(\int_E|f|^p)^{1/p}$ stands for $\|f\chi_E\|_\infty$.
	
	The Hardy--Littlewood maximal operator is
	$Mf(x)=\sup_{B\ni x}\frac1{|B|}\int_B|f|$. For $1<p<\infty$, the
	Muckenhoupt class $\mathcal A_p$ consists of the weights $w$ such that
	\[
	[w]_{\mathcal{A}_p}:=\sup_{B}
	\Bigl(\frac{1}{|B|}\int_B w\Bigr)
	\Bigl(\frac{1}{|B|}\int_B w^{-\frac{1}{p-1}}\Bigr)^{p-1}<\infty .
	\]
	A weight $w$ belongs to $\mathcal A_1$ if $Mw\le c\,w$ a.e.;
	equivalently, there exists $C>0$ such that, for every ball $B$,
	\[
	\frac{1}{|B|}\int_B w\le C\,\einf_{B} w .
	\]
	We set $\mathcal A_\infty=\bigcup_{p\ge1}\mathcal A_p$ and denote by
	$\tilde q_w=\inf\{q>1: w\in\mathcal A_q\}$ the critical index of $w$.
	
	A weight $w$ satisfies the reverse H\"older inequality with exponent
	$s>1$, $w\in RH_s$, if there is $C>0$ such that for every ball $B$
	\[
	\Bigl(\frac{1}{|B|}\int_B w^{s}\Bigr)^{1/s}\le\frac{C}{|B|}\int_B w .
	\]
	If $w\in RH_s$ then $w\in RH_{s'}$ for $1<s'<s$ and $w\in RH_{s+\epsilon}$
	for some $\epsilon>0$. We denote by $r_w=\sup\{r>1:w\in RH_r\}$ the
	critical index for the reverse H\"older condition. It is well known that
	$w\in\mathcal A_\infty$ if and only if $w\in RH_s$ for some $s>1$; thus
	$1<r_w\le\infty$ for all $w\in\mathcal A_\infty$.
	
	We shall use the following well-known properties of $\mathcal A_1$ weights;
	see, for instance, \cite{GR}, \cite{D} or \cite{G}.
	
	\begin{lemma}\label{L-A1}
		Let $v\in\mathcal A_1$.
		\begin{itemize}
			\item[(i)] $v^{\theta}\in\mathcal A_1$ for every $0<\theta\le1$.
			\item[(ii)] There exists $\varepsilon>0$ such that $v^{1+\varepsilon}\in\mathcal A_1$.
			\item[(iii)] $v^{\theta}\in RH_{1/\theta}$ for every $0<\theta<1$.
			\item[(iv)] $v(\lambda B)\le C\lambda^n v(B)$ for every ball $B$ and
			$\lambda\ge1$.
		\end{itemize}
	\end{lemma}
	
	Property (i) is a consequence of Jensen's inequality, (ii) follows from the
	reverse H\"older inequality satisfied by $\mathcal A_\infty$ weights, and
	(iii) and (iv) follow directly from the definition of $\mathcal A_1$.
	
	\subsection{Weighted Hardy spaces}
	Topologize $\mathcal{S}(\RR^n)$ by the semi-norms
	\[
	\|\phi\|_{\mu,\kappa}=\sup_{x}|x^\mu\partial^\kappa\phi(x)|.
	\]
	For $N\in\mathbb N$ let
	\[
	S_N=\{\phi\in\mathcal{S}(\RR^n):\|\phi\|_{\mu,\kappa}\le1,\
	|\mu|,|\kappa|\le N\},
	\]
	and for $f\in\mathcal S'(\RR^n)$ let
	\[
	M_Nf(x)=\sup_{t>0}\sup_{\phi\in S_N}\bigl|t^{-n}\phi(t^{-1}\cdot)*f(x)\bigr|.
	\]
	Given $w\in\mathcal A_\infty$ and $p>0$, the weighted Hardy space
	$H_w^p(\RR^n)$ consists of all $f\in\mathcal S'(\RR^n)$ such that
	$\|f\|_{H_w^p}=\|M_Nf\|_{L_w^p}<\infty$, for $N$ sufficiently large
	(see \cite[Chapter~VI]{ST}). Let
	\[
	\mathcal D_0=\{\phi\in\mathcal S(\RR^n):\widehat\phi\in C_c^\infty(\RR^n),\
	0\notin\operatorname{supp}\widehat\phi\}
	\]
	(see \cite[Chapter~V]{ST}).
	
	\begin{theorem}[{\cite[Chapter~VII, Theorem~1]{ST}}]\label{T-dense}
		Let $w$ be a doubling weight on $\RR^{n}$. Then $\mathcal D_0$ is dense in
		$H^{p}_w(\RR^{n})$, $0<p<\infty$.
	\end{theorem}
	
	We recall the definition of $w$-$(p,p_0,d)$ atom given in \cite{Rocha}.
	Let $w\in\mathcal A_\infty$, $0<p\le1$,
	$\max\{1,p\frac{r_w}{r_w-1}\}<p_0\le\infty$ (where
	$\frac{r_w}{r_w-1}=1$ if $r_w=\infty$), and $d\in\mathbb Z$ with
	$d\ge\lfloor n(\frac{\tilde q_w}{p}-1)\rfloor$. A function $a$ is a
	$w$-$(p,p_0,d)$ atom centered at $x_0$ if
	\begin{itemize}
		\item[(a1)] $a\in L^{p_0}(\RR^{n})$ is supported in a ball $B=B(x_0,r)$;
		\item[(a2)] $\|a\|_{L^{p_0}}\le|B|^{1/p_0}w(B)^{-1/p}$;
		\item[(a3)] $\int x^{\mu}a(x)\,dx=0$ for every multi-index $\mu$ with
		$|\mu|\le d$.
	\end{itemize}
	
	\begin{lemma}[{\cite[Lemma~2.8]{Rocha}}]\label{Lemma8}
		There exists $C>0$ such that $\|a\|_{H^{p}_w}\le C$ for every
		$w$-$(p,p_0,d)$ atom $a$.
	\end{lemma}
	
	\begin{theorem}[{\cite[Theorem~2.9]{Rocha}}]\label{Theorem9}
		Let $w\in\mathcal A_\infty$, $0<p\le1$, let $p_0$ and $d$ be as in the
		definition of $w$-$(p,p_0,d)$ atom, and let $f\in\mathcal D_0$. There exist
		$w$-$(p,p_0,d)$ atoms $a_j$ and scalars $\lambda_j$ with
		$\sum_j|\lambda_j|^{p}\le c\|f\|^{p}_{H^{p}_w}$ such that
		$f=\sum_j\lambda_ja_j$, where the series converges in $L^{s}(\RR^{n})$ for
		every $1<s<\infty$ and pointwise.
	\end{theorem}
	
	The atoms constructed in the proof of \cite[Theorem~2.9]{Rocha} satisfy
	$\|a\|_{L^\infty}\le w(B)^{-1/p}$, so they are $w$-$(p,p_0,d)$ atoms for
	every admissible $p_0$; moreover, the order $d$ of their vanishing moments can
	be prescribed.
	
	\section{Size and strengthened H\"ormander conditions}\label{S-Hormander}
	
	\subsection{Definitions and basic properties}
	The normalized annular norms introduced in Section~1 are invariant under
	dilations: if $f_\lambda(x)=f(x/\lambda)$, then
	$\|f_\lambda\|_{\tau,|x|\sim\lambda s}=\|f\|_{\tau,|x|\sim s}$. By
	H\"older's inequality, $\|f\|_{\tau_1,|x|\sim s}\le\|f\|_{\tau_2,|x|\sim s}$
	whenever $1\le\tau_1\le\tau_2\le\infty$.
	
	\begin{definition}[Size condition, {\cite[Definition~3.1]{BLRi}}]\label{def-S}
		Let $1\le\tau\le\infty$ and $0\le\beta<n$. A measurable function $k$
		belongs to $S_{\beta,\tau}$ if there exists $C>0$ such that
		\[
		\|k\|_{\tau,|x|\sim s}\le C\,s^{\beta-n}\qquad\text{for all } s>0.
		\]
	\end{definition}
	
	We now recall the strengthened fractional H\"ormander condition
	considered in this paper, which was introduced by Han and Wu
	\cite[Definition~5.6]{HW} under the name of fractional $L^\tau$-H\"ormander
	condition of order $\delta$. It is obtained from the fractional H\"ormander
	condition of \cite{BLRi} by adding the factor $2^{j\delta}$, which is adapted
	to the $H^p$--$L^q$ setting of Theorem~\ref{T-HpLq}: for
	$\delta=n(\frac1p-1)$, it is exactly the factor that appears when $T_\alpha$
	is applied to the atoms of $H^p_{w^p}$ (see Remark~\ref{R-origin}).
	
	\begin{definition}[{\cite[Definition~5.6]{HW}}]\label{def-Hdelta}
		Let $1\le\tau\le\infty$, $0\le\beta<n$, and $\delta\ge0$. A measurable
		function $k$ belongs to $H_{\beta,\tau}^{(\delta)}$ if there exist $c\ge1$ and
		$C>0$ such that, for every $h\in\RR^n$ and $R>c|h|$,
		\[
		\sum_{j=1}^\infty 2^{j\delta}(2^jR)^{n-\beta}
		\|k(\cdot-h)-k\|_{\tau,|x|\sim2^jR}\le C.
		\]
		For $\delta=0$ this is the fractional H\"ormander condition $H_{\beta,\tau}$
		of \cite[Definition~3.2]{BLRi} (the case $k=0$ there).
	\end{definition}
	
	The formulation in \cite{HW} differs from Definition~\ref{def-Hdelta} only in
	the normalization of the averages and in the constant $c$, and there the
	order is restricted to $0\le\delta\le1$.
	
	With this notation, $(\mathrm H_1)$ and $(\mathrm H_2)$ say that
	$k_i\in S_{n-n/q_i,\,p_i}$ and $k_i\in H^{(\delta)}_{n-n/q_i,\,p_i}$ for
	$1\le i\le m$.
	
	\begin{remark}\label{R-inclusions}
		Since $2^{j\delta_1}\le 2^{j\delta_2}$ for $j\ge1$ and by the monotonicity
		of the normalized norms in $\tau$, we have
		\[
		H_{\beta,\tau_2}^{(\delta_2)}\subseteq H_{\beta,\tau_1}^{(\delta_1)}
		\subseteq H_{\beta,\tau_1},
		\qquad S_{\beta,\tau_2}\subseteq S_{\beta,\tau_1},
		\]
		whenever $1\le\tau_1\le\tau_2\le\infty$ and $0\le\delta_1\le\delta_2$.
	\end{remark}
	
	\begin{example}[Homogeneous kernels]\label{Ex-homog}
		Let $\Omega_i\in L^{p_i}(S^{n-1})$, $q_i<p_i\le\infty$, be homogeneous of
		degree zero and $k_i(x)=\Omega_i(x)|x|^{-n/q_i}$. If $p_i<\infty$, polar
		coordinates give
		\[
		\|k_i\|_{p_i,|x|\sim s}^{p_i}
		=\frac{\|\Omega_i\|^{p_i}_{L^{p_i}(S^{n-1})}}{|B(0,2s)|}
		\int_s^{2s}r^{n-1-np_i/q_i}\,dr
		=C\,s^{-np_i/q_i},
		\]
		and the case $p_i=\infty$ is immediate. Hence $k_i\in S_{n-n/q_i,p_i}$ and
		$(\mathrm H_1)$ holds. This is the setting of \cite{RiU2}.
	\end{example}
	
	The following consequence of the size condition replaces the polar
	coordinates computation available for homogeneous kernels.
	
	\begin{lemma}\label{L-tail}
		Let $1<q_i<p_i\le\infty$, $\gamma_i=\frac n{q_i}-\frac n{p_i}>0$, and
		$k_i\in S_{n-n/q_i,p_i}$. Then there is $C>0$ such that
		\[
		\|k_i\chi_{\{|x|>u\}}\|_{L^{p_i}}\le C\,u^{-\gamma_i}
		\qquad\text{for all } u>0 .
		\]
	\end{lemma}
	
	\begin{proof}
		For $l\ge0$, by the size condition $k_i\in S_{n-n/q_i,p_i}$,
		\[
		\|k_i\chi_{\{2^lu<|x|\le2^{l+1}u\}}\|_{L^{p_i}}
		=|B(0,2^{l+1}u)|^{1/p_i}\,\|k_i\|_{p_i,|x|\sim2^lu}
		\le C(2^lu)^{\frac n{p_i}-\frac n{q_i}}=C(2^lu)^{-\gamma_i},
		\]
		with the obvious interpretation if $p_i=\infty$. Summing over $l\ge0$ and
		using the triangle inequality in $L^{p_i}$ gives the claim, since
		$\gamma_i>0$.
	\end{proof}
	
	\begin{remark}\label{R-delta-large}
		The restriction $\delta<1$ is also natural from the point of view of the
		hypotheses: if $\delta>1$, then
		$S_{\beta,\tau}\cap H^{(\delta)}_{\beta,\tau}=\{0\}$ for every
		$1\le\tau\le\infty$ and $0\le\beta<n$, so that $(\mathrm H_1)$ and
		$(\mathrm H_2)$ would force $T_\alpha\equiv0$. Indeed, let
		$k\in S_{\beta,\tau}\cap H^{(\delta)}_{\beta,\tau}$ with constants $c$ and $C$,
		fix $s>0$ and let $0<|h|<s/(2c)$. Choose $j\ge1$ with
		$s/(c2^{j+1})\le|h|<s/(c2^j)$ and apply Definition~\ref{def-Hdelta} with
		$R=2^{-j}s>c|h|$. Keeping only the $j$-th term of the sum, we get
		\[
		\|k(\cdot-h)-k\|_{\tau,|x|\sim s}\le Cs^{\beta-n}2^{-j\delta}
		\le Cs^{\beta-n}\Bigl(\frac{2c|h|}{s}\Bigr)^{\delta},
		\]
		which is $o(|h|)$ as $h\to0$, since $\delta>1$. Hence, if
		$\varphi\in C_c^\infty(\{s<|x|<2s\})$ and $1\le l\le n$, H\"older's
		inequality gives
		\[
		\Bigl|\int k\,\partial_l\varphi\Bigr|
		=\lim_{t\to0}\Bigl|\int\frac{k(x-te_l)-k(x)}{t}\,\varphi(x)\,dx\Bigr|
		\le C_{s,\varphi}\lim_{t\to0}
		\frac{\|k(\cdot-te_l)-k\|_{\tau,|x|\sim s}}{|t|}=0 .
		\]
		Since the open annuli $\{s<|x|<2s\}$, $s>0$, cover $\RR^n\setminus\{0\}$, we
		get $\nabla k=0$ in $\mathcal D'(\RR^n\setminus\{0\})$, so $k$ coincides a.e.
		with a constant $c_V$ on each connected component $V$ of
		$\RR^n\setminus\{0\}$. As $|V\cap\{|x|\sim s\}|/|B(0,2s)|$ does not depend on
		$s$, the size condition gives $|c_V|\le Cs^{\beta-n}$ for every $s>0$, and
		letting $s\to\infty$ we obtain $c_V=0$. In particular, the restriction
		$\delta\le1$ in \cite[Definition~5.6]{HW} entails no loss of generality.
	\end{remark}
	
	By Remark~\ref{R-inclusions}, the classes $H^{(\delta)}_{\beta,\tau}$
	decrease as $\delta$ increases. The next result shows that these inclusions
	are strict. In part~(a) we consider the homogeneous kernels
	$k_\varepsilon(x)=\Omega_\varepsilon(x)|x|^{\beta-n}$, with
	$\Omega_\varepsilon(x)=|x_n/|x||^{\varepsilon}$, and determine exactly the
	values of $\delta$ for which $k_\varepsilon\in H^{(\delta)}_{\beta,\tau}$:
	these are the $\delta<\varepsilon+\frac1\tau$. For $\varepsilon>0$, the
	function $\Omega_\varepsilon$ is H\"older continuous of order $\varepsilon$
	and vanishes on the hyperplane $\{x_n=0\}$, while for $\varepsilon<0$ it is
	singular there; the term $\frac1\tau$ appears because the irregularity of
	$k_\varepsilon$ is concentrated near this hyperplane. Since the threshold
	$\varepsilon+\frac1\tau$ takes every value in $(0,1)$ except $\frac1\tau$,
	any two classes $H^{(\delta_1)}_{\beta,\tau}$ and
	$H^{(\delta_2)}_{\beta,\tau}$ with $\delta_1<\delta_2$, $\delta_1<1$, are
	separated by some $k_\varepsilon$; this is part~(b).
	
	\begin{proposition}\label{P-strict}
		Let $n\ge2$, $0\le\beta<n$ and $1\le\tau\le\infty$, where $\frac1\tau=0$ if
		$\tau=\infty$.
		\begin{itemize}
			\item[(a)] Let $-\frac1\tau<\varepsilon<1-\frac1\tau$, $\varepsilon\neq0$, and
			$k_\varepsilon(x)=|x_n|^{\varepsilon}|x|^{\beta-n-\varepsilon}$. Then
			$k_\varepsilon\in S_{\beta,\tau}$, and $k_\varepsilon\in H^{(\delta)}_{\beta,\tau}$
			if and only if $\delta<\varepsilon+\frac1\tau$.
			\item[(b)] If $0\le\delta_1<\delta_2$ and $\delta_1<1$, then
			\[
			S_{\beta,\tau}\cap H^{(\delta_2)}_{\beta,\tau}\subsetneq
			S_{\beta,\tau}\cap H^{(\delta_1)}_{\beta,\tau}.
			\]
		\end{itemize}
	\end{proposition}
	
	\begin{proof}
		(a) We write $x=(\tilde x,x_n)\in\RR^{n-1}\times\RR$. Let
		$\theta=\varepsilon+\frac1\tau$, so that $0<\theta<1$, and write
		$k_\varepsilon(x)=g(x_n)\varrho(x)$, with $g(u)=|u|^{\varepsilon}$ and
		$\varrho(x)=|x|^{\beta-n-\varepsilon}$. Since $\varepsilon\tau>-1$, the function
		$|x_n/|x||^{\varepsilon}$ belongs to $L^\tau(S^{n-1})$, and
		$k_\varepsilon\in S_{\beta,\tau}$ as in Example~\ref{Ex-homog}.
		
		We use two elementary facts. First, for $b\in\RR$,
		\begin{equation}\label{strict-1d}
			\Bigl(\int_{\RR}\bigl||u-b|^{\varepsilon}-|u|^{\varepsilon}\bigr|^{\tau}du\Bigr)^{1/\tau}
			=C_{\varepsilon,\tau}|b|^{\varepsilon+\frac1\tau},
			\qquad C_{\varepsilon,\tau}<\infty .
		\end{equation}
		Indeed, the change of variables $u=|b|v$ reduces \eqref{strict-1d} to
		$b=\pm1$, and the integral converges because $\varepsilon\tau>-1$ (near
		$v=0$ and $v=\pm1$) and $(\varepsilon-1)\tau<-1$ (at infinity, where the
		integrand is $O(|v|^{(\varepsilon-1)\tau})$); if $\tau=\infty$, then
		$0<\varepsilon<1$ and
		$\bigl||u-b|^{\varepsilon}-|u|^{\varepsilon}\bigr|\le|b|^{\varepsilon}$.
		Second, since $\{|x|\le2s\}\subset\{|\tilde x|\le2s,\ |x_n|\le2s\}$, for every
		function $F$ on $\RR$ we have
		\begin{equation}\label{strict-xn}
			\|F(x_n)\|_{\tau,|x|\sim s}
			\le C\Bigl(\frac1s\int_{-2s}^{2s}|F(u)|^\tau\,du\Bigr)^{1/\tau}.
		\end{equation}
		
		\smallskip
		We first show that $k_\varepsilon\in H^{(\delta)}_{\beta,\tau}$ if
		$\delta<\theta$. Let $h\in\RR^n$,
		$t=|h|$, $R>4t$, $j\ge1$ and $s=2^jR$, so that $t<s/8$. If $|x|\sim s$, every
		point of the segment joining $x-h$ and $x$ has modulus between $7s/8$ and
		$3s$, so $|\varrho(x-h)|\le Cs^{\beta-n-\varepsilon}$ and, by the mean value
		theorem, $|\varrho(x-h)-\varrho(x)|\le C\frac ts\,s^{\beta-n-\varepsilon}$.
		Writing
		\[
		k_\varepsilon(x-h)-k_\varepsilon(x)
		=[g(x_n-h_n)-g(x_n)]\varrho(x-h)+g(x_n)[\varrho(x-h)-\varrho(x)],
		\]
		and using \eqref{strict-xn}, \eqref{strict-1d} with $b=h_n$, and
		$\int_{-2s}^{2s}|u|^{\varepsilon\tau}du=Cs^{\varepsilon\tau+1}$, we get
		\[
		s^{n-\beta}\|k_\varepsilon(\cdot-h)-k_\varepsilon\|_{\tau,|x|\sim s}
		\le C\Bigl(\frac ts\Bigr)^{\theta}+C\frac ts
		\le C\Bigl(\frac ts\Bigr)^{\theta},
		\]
		since $\theta<1$ and $t<s$. Hence
		\[
		\sum_{j=1}^\infty2^{j\delta}(2^jR)^{n-\beta}
		\|k_\varepsilon(\cdot-h)-k_\varepsilon\|_{\tau,|x|\sim2^jR}
		\le C\Bigl(\frac tR\Bigr)^{\theta}\sum_{j=1}^\infty2^{j(\delta-\theta)}\le C,
		\]
		so $k_\varepsilon\in H^{(\delta)}_{\beta,\tau}$, with $c=4$.
		
		\smallskip
		We now show that $k_\varepsilon\notin H^{(\delta)}_{\beta,\tau}$ if
		$\delta\ge\theta$. Suppose
		that $k_\varepsilon\in H^{(\delta)}_{\beta,\tau}$ with constants $c$ and $C$.
		Let $t>0$, $h=te_n$, $R>\max\{c,4\}t$, $j\ge1$ and $s=2^jR$. The set
		\[
		E_s=\{x:\ 6s/5<|\tilde x|<9s/5,\ |x_n|<t/8\}
		\]
		is contained in $\{s<|x|\le2s\}$ and, since $n\ge2$, $|E_s|\ge c_ns^{n-1}t$.
		For $x\in E_s$ we have $|x_n|<t/8$ and $7t/8<|x_n-t|<9t/8$, hence, as
		$\varepsilon\ne0$ and $|\varepsilon|<1$,
		\[
		|g(x_n-t)-g(x_n)|\ge c_\varepsilon t^{\varepsilon},\qquad
		g(x_n-t)\le2t^{\varepsilon},\qquad
		c_\varepsilon=\Bigl|\Bigl(\frac78\Bigr)^{\varepsilon}-\Bigl(\frac18\Bigr)^{\varepsilon}\Bigr|>0 .
		\]
		Together with $|\varrho(x-h)-\varrho(x)|\le C\frac ts\,\varrho(x)$, this gives,
		for $x\in E_s$,
		\[
		|k_\varepsilon(x-h)-k_\varepsilon(x)|
		\ge t^{\varepsilon}\varrho(x)\Bigl(c_\varepsilon-C\frac ts\Bigr)
		\ge\frac{c_\varepsilon}2\,t^{\varepsilon}\varrho(x)
		\ge c\,t^{\varepsilon}s^{\beta-n-\varepsilon}
		\]
		for every $j\ge j_0$, where $j_0$ depends only on $n$, $\beta$ and
		$\varepsilon$, since $t/s<2^{-j}$. As $|E_s|/|B(0,2s)|\ge c\,t/s$, we
		obtain, for $j\ge j_0$,
		\[
		2^{j\delta}(2^jR)^{n-\beta}\|k_\varepsilon(\cdot-h)-k_\varepsilon\|_{\tau,|x|\sim2^jR}
		\ge c\,2^{j\delta}\Bigl(\frac t{2^jR}\Bigr)^{\theta}
		=c\Bigl(\frac tR\Bigr)^{\theta}2^{j(\delta-\theta)}
		\ge c\Bigl(\frac tR\Bigr)^{\theta},
		\]
		so the series in Definition~\ref{def-Hdelta} diverges, a contradiction.
		This proves (a).
		
		\smallskip\noindent(b) Choose $\theta$ with
		$\delta_1<\theta<\min\{\delta_2,1\}$ and $\theta\neq\frac1\tau$, and let
		$\varepsilon=\theta-\frac1\tau$. Then $-\frac1\tau<\varepsilon<1-\frac1\tau$
		and $\varepsilon\neq0$, so, by (a), $k_\varepsilon$ belongs to
		$S_{\beta,\tau}\cap H^{(\delta_1)}_{\beta,\tau}$ but not to
		$H^{(\delta_2)}_{\beta,\tau}$, since $\delta_2>\theta$.
	\end{proof}
	
	\subsection{Sufficient conditions for \texorpdfstring{$(\mathrm H_2)$}{(H2)}}
	The following results show that $(\mathrm H_2)$ is satisfied by large
	classes of kernels, homogeneous or not. In all of them the restriction
	$\delta<1$ of Theorem~\ref{T-HpLq} plays a role.
	
	\begin{proposition}[Pointwise smoothness]\label{P-Holder}
		Let $0\le\beta<n$, $0<\varepsilon\le1$, and let $k$ satisfy, for some
		$c\ge1$,
		\[
		|k(x-h)-k(x)|\le C\frac{|h|^{\varepsilon}}{|x|^{n-\beta+\varepsilon}},
		\qquad |x|>c|h| .
		\]
		Then $k\in H^{(\delta)}_{\beta,\infty}$ for every $0\le\delta<\varepsilon$.
		In particular, if $\varepsilon=1$ (the condition $H^*_{\beta,\infty}$ of
		\cite[Definition~3.3]{BLRi}), then $k\in H^{(\delta)}_{\beta,\infty}$ for
		every $0\le\delta<1$.
	\end{proposition}
	
	\begin{proof}
		Let $R>c|h|$ and $j\ge1$. If $|x|>2^jR$, then $|x|>c|h|$, and hence
		$\|k(\cdot-h)-k\|_{\infty,|x|\sim2^jR}\le C|h|^\varepsilon(2^jR)^{\beta-n-\varepsilon}$.
		Therefore
		\[
		\sum_{j=1}^\infty2^{j\delta}(2^jR)^{n-\beta}\|k(\cdot-h)-k\|_{\infty,|x|\sim2^jR}
		\le C\Bigl(\frac{|h|}{R}\Bigr)^{\varepsilon}\sum_{j=1}^\infty2^{j(\delta-\varepsilon)}
		\le C,
		\]
		since $|h|<R$ and $\delta<\varepsilon$.
	\end{proof}
	
	The next result is the analogue of \cite[Proposition~4.1]{BLRi} for the
	classes $H^{(\delta)}_{\beta,\tau}$. It allows us to produce fractional
	kernels from kernels of singular integral type.
	
	\begin{proposition}\label{P-BLRi35}
		Let $1\le\tau\le\infty$, $0<\beta<n$ and $0\le\delta<1$. If
		$K\in S_{0,\tau}\cap H^{(\delta)}_{0,\tau}$, then
		$k(x)=|x|^{\beta}K(x)$ belongs to $S_{\beta,\tau}\cap H^{(\delta)}_{\beta,\tau}$.
	\end{proposition}
	
	\begin{proof}
		Since $|x|^\beta\le(2s)^\beta$ on $\{|x|\sim s\}$, we have
		$\|k\|_{\tau,|x|\sim s}\le Cs^{\beta}\|K\|_{\tau,|x|\sim s}\le Cs^{\beta-n}$,
		so $k\in S_{\beta,\tau}$.
		
		Let $c\ge1$ be the constant of $K$, $R>c|h|$, $j\ge1$ and $s=2^jR$. If
		$|x|\sim s$, then $|h|<R\le s/2$, so $s/2<|x-h|<5s/2$ and, by the mean
		value theorem, $\bigl||x-h|^\beta-|x|^\beta\bigr|\le Cs^{\beta-1}|h|$. Hence
		\begin{align*}
			|k(x-h)-k(x)|
			&\le|x-h|^\beta|K(x-h)-K(x)|+|K(x)|\,\bigl||x-h|^\beta-|x|^\beta\bigr|\\
			&\le Cs^\beta|K(x-h)-K(x)|+Cs^{\beta}\frac{|h|}{s}|K(x)|,
		\end{align*}
		and, using $K\in S_{0,\tau}$,
		\[
		\|k(\cdot-h)-k\|_{\tau,|x|\sim s}
		\le Cs^\beta\|K(\cdot-h)-K\|_{\tau,|x|\sim s}+Cs^{\beta-n}\frac{|h|}{s}.
		\]
		Multiplying by $2^{j\delta}s^{n-\beta}$ and summing over $j\ge1$,
		\begin{align*}
			&\sum_{j=1}^\infty2^{j\delta}(2^jR)^{n-\beta}\|k(\cdot-h)-k\|_{\tau,|x|\sim2^jR}\\
			&\qquad\le C\sum_{j=1}^\infty2^{j\delta}(2^jR)^{n}\|K(\cdot-h)-K\|_{\tau,|x|\sim2^jR}
			+C\frac{|h|}{R}\sum_{j=1}^\infty2^{j(\delta-1)}\le C,
		\end{align*}
		since $K\in H^{(\delta)}_{0,\tau}$, $|h|<R$ and $\delta<1$.
	\end{proof}
	
	For rough homogeneous kernels, $(\mathrm H_2)$ follows from a Dini type
	condition of order $\delta$. For $\Omega\in L^\tau(S^{n-1})$ let
	\[
	\varpi_\tau(u)=\sup_{|\varrho|\le u}
	\|\Omega(\varrho\,\cdot)-\Omega\|_{L^\tau(S^{n-1})},
	\]
	where the supremum is taken over rotations $\varrho$ with
	$|\varrho|=\|\varrho-I\|\le u$.
	
	\begin{proposition}\label{P-Dini}
		Let $1<\tau<\infty$, $0<\beta<n$, $0\le\delta<1$, and let
		$\Omega\in L^\tau(S^{n-1})$ be homogeneous of degree zero with
		\begin{equation}\label{Dini-delta}
			\int_0^1\varpi_\tau(u)\,\frac{du}{u^{1+\delta}}<\infty .
		\end{equation}
		Then $k(x)=\Omega(x)|x|^{\beta-n}$ belongs to
		$S_{\beta,\tau}\cap H^{(\delta)}_{\beta,\tau}$.
	\end{proposition}
	
	\begin{proof}
		That $k\in S_{\beta,\tau}$ follows as in Example~\ref{Ex-homog}. For the
		smoothness we use \cite[Lemma 4]{DL2} (whose proof follows
		\cite[Lemma 5]{KuW}): since \eqref{Dini-delta} implies the $L^\tau$-Dini
		condition $\int_0^1\varpi_\tau(u)\,\frac{du}{u}<\infty$, and after dividing
		by $|B(0,2s)|^{1/\tau}$, there is $C>0$ such that, for $s>4|h|$,
		\begin{equation}\label{DL-lemma}
			\|k(\cdot-h)-k\|_{\tau,|x|\sim s}
			\le Cs^{\beta-n}\Bigl(\frac{|h|}{s}
			+\int_{|h|/2s}^{|h|/s}\varpi_\tau(u)\,\frac{du}{u}\Bigr).
		\end{equation}
		
		Let $R>2|h|$, $j\ge1$ and $s=2^jR$, so that $s>4|h|$. If
		$u\in[|h|/2^{j+1}R,\,|h|/2^jR]$, then $2^{j}\le|h|/(Ru)\le u^{-1}$, so
		$2^{j\delta}\le u^{-\delta}$. Hence, by \eqref{DL-lemma},
		\[
		\sum_{j=1}^\infty2^{j\delta}(2^jR)^{n-\beta}\|k(\cdot-h)-k\|_{\tau,|x|\sim2^jR}
		\le C\frac{|h|}{R}\sum_{j=1}^\infty2^{j(\delta-1)}
		+C\int_0^{1}\varpi_\tau(u)\,\frac{du}{u^{1+\delta}}\le C .
		\]
	\end{proof}
	
	For $\delta=0$, Proposition~\ref{P-Dini} is essentially
	\cite[Proposition~4.2]{BLRi}, with $k=0$ there and the Young function
	$t\mapsto t^\tau$. For $0\le\delta<1$, the conclusion of
	Proposition~\ref{P-Dini} is stated without proof in \cite[Remark~5.7]{HW}.
	
	\begin{remark}
		For $m=1$, Ding and Lu \cite[Theorem 2]{DL2} proved the (unweighted)
		$H^p-L^q$ boundedness of $T_{\Omega,\alpha}$, for $0<\alpha<1$ and
		$\frac n{n+\alpha}\le p<1$, for $\Omega\in L^r(S^{n-1})$ with
		$r>\frac n{n-\alpha}$ satisfying the Dini condition
		$\int_0^1\varpi_r(u)\,u^{-1-\alpha}\,du<\infty$. Since
		$\delta=n(\frac1p-1)\le\alpha$ in that range, condition \eqref{Dini-delta},
		of order $\delta$, is weaker than theirs; in fact, their argument (see
		\cite[(3.6)]{DL2}) only uses the Dini condition of order $\delta$. Thus, for
		$m=1$, Theorem~\ref{T-HpLq} and Proposition~\ref{P-Dini} give a weighted
		version of \cite[Theorem 2]{DL2} for $\frac n{n+\alpha}<p<1$, under the Dini
		condition of order $\delta$.
	\end{remark}
	
	\begin{remark}
		If $\Omega$ is H\"older continuous of order $\varepsilon$ on $S^{n-1}$, then
		$\varpi_\tau(u)\le Cu^{\varepsilon}$ for every $\tau$, and
		\eqref{Dini-delta} holds whenever $\varepsilon>\delta$; in this case
		Proposition~\ref{P-Holder} also gives $k\in H^{(\delta)}_{\beta,\infty}$. Compare with Remark~\ref{R-counterexample},
		where $\Omega_1$ is H\"older of order $\varepsilon\le\delta-1/q$ and the
		conclusion of Theorem~\ref{T-HpLq} fails.
	\end{remark}
	
	\begin{example}[A non-homogeneous kernel]\label{Ex-nonhomog}
		Let $\Omega\in C^1(S^{n-1})$, extended homogeneously of degree zero, and
		\[
		k_i(x)=\frac{\Omega(x)\,\cos(\log|x|)}{|x|^{n/q_i}} .
		\]
		Then $|k_i(x)|\le C|x|^{-n/q_i}$ and $|\nabla k_i(x)|\le C|x|^{-n/q_i-1}$
		for $x\neq0$. If $|x|>2|h|$, the segment from $x-h$ to $x$ stays in
		$\{|z|>|x|/2\}$, so $|k_i(x-h)-k_i(x)|\le C|h|\,|x|^{-n/q_i-1}$. Hence
		$k_i\in S_{n-n/q_i,\infty}$ and, by Proposition~\ref{P-Holder},
		$k_i\in H^{(\delta)}_{n-n/q_i,\infty}$ for every $0\le\delta<1$. Thus
		$(\mathrm H_1)$ and $(\mathrm H_2)$ hold with $p_i=\infty$, although $k_i$ is
		not homogeneous. More generally, by Proposition~\ref{P-BLRi35}, one may take
		$k_i(x)=|x|^{n-n/q_i}K_i(x)$ with $K_i\in S_{0,\infty}\cap H^{(\delta)}_{0,\infty}$,
		for instance a Calder\'on--Zygmund kernel with H\"older regularity of order
		$\varepsilon>\delta$. Other non-homogeneous kernels satisfying size and
		H\"ormander conditions of Young type are given, in dimension one, in
		\cite[Section 4]{FR}.
	\end{example}
	
	\subsection{The homogeneous case}\label{S-homog}
	In this subsection $k_i(x)=\Omega_i(x)|x|^{-n/q_i}$ as in
	Example~\ref{Ex-homog}. Definition~\ref{def-Hdelta} makes sense for every
	$\beta\in\RR$, and here we allow any real $\beta$. For $\beta,\delta\in\RR$,
	$h\in\RR^n$ and $R>0$ set
	\[
	\Lambda^{(\delta)}_\beta(h,R)=\sum_{j=1}^\infty 2^{j\delta}(2^jR)^{n-\beta}
	\|k_i(\cdot-h)-k_i\|_{p_i,|x|\sim2^jR},
	\qquad \Lambda_\beta=\Lambda_\beta^{(0)} .
	\]
	The next result shows that, among classical fractional H\"ormander classes,
	the only one compatible with the homogeneity of $k_i$ is
	$H_{n-n/q_i,p_i}$. Consequently, the extra decay required in
	$(\mathrm H_2)$ cannot be obtained by changing the index $\beta$, and has to
	be imposed through the factor $2^{j\delta}$.
	
	\begin{theorem}\label{T-index}
		Let $1\le i\le m$, $1\le p_i\le\infty$, and $k_i(x)=\Omega_i(x)|x|^{-n/q_i}$
		with $\Omega_i\in L^1(S^{n-1})$ homogeneous of degree zero and
		$k_i\not\equiv0$.
		\begin{itemize}
			\item[(a)] If $k_i\in H_{\beta,p_i}$ for some $\beta\in\RR$, then
			$\beta=n-n/q_i$.
			\item[(b)] In particular, $k_i\notin H_{n-n/q_i-\delta,\,p_i}$ for every
			$\delta\neq0$.
			\item[(c)] For every $\delta\in\RR$, $\Lambda^{(\delta)}_{n-n/q_i}(\lambda h,\lambda R)
			=\Lambda^{(\delta)}_{n-n/q_i}(h,R)$ for all $\lambda>0$.
		\end{itemize}
	\end{theorem}
	
	\begin{proof}
		The proof relies on the behavior of the sums $\Lambda^{(\delta)}_\beta$ under
		dilations, which is a consequence of the homogeneity of $k_i$. For $h\in\RR^n$,
		$R>0$ and $j\ge1$, let $N_j(h,R)=\|k_i(\cdot-h)-k_i\|_{p_i,|x|\sim2^jR}$.
		Since $k_i$ is homogeneous of degree $-n/q_i$, we have
		$k_i(\lambda x-\lambda h)-k_i(\lambda x)=\lambda^{-n/q_i}(k_i(x-h)-k_i(x))$
		for a.e. $x$ and every $\lambda>0$. If $p_i<\infty$, the change of variables
		$x=\lambda y$ and the identity $|B(0,2^{j+1}\lambda R)|=\lambda^n|B(0,2^{j+1}R)|$
		give
		\begin{align*}
			N_j(\lambda h,\lambda R)
			&=\Bigl(\frac{\lambda^n\lambda^{-np_i/q_i}}{\lambda^n|B(0,2^{j+1}R)|}
			\int_{2^jR<|y|\le2^{j+1}R}|k_i(y-h)-k_i(y)|^{p_i}\,dy\Bigr)^{1/p_i}\\
			&=\lambda^{-n/q_i}N_j(h,R),
		\end{align*}
		and the same identity holds for $p_i=\infty$, taking essential suprema.
		Multiplying by $2^{j\delta}(2^j\lambda R)^{n-\beta}$ and summing in $j$, we
		obtain
		\begin{equation}\label{T-scaling}
			\Lambda^{(\delta)}_\beta(\lambda h,\lambda R)=\lambda^{\theta}\Lambda^{(\delta)}_\beta(h,R),
			\qquad\theta=n-\beta-\frac n{q_i}.
		\end{equation}
		When $\beta=n-n/q_i$ we have $\theta=0$, and \eqref{T-scaling} is precisely
		the statement of (c).
		
		We now prove (a). Assume that $k_i\in H_{\beta,p_i}$ with constants $c$ and
		$C$ and suppose, by contradiction, that $\beta\ne n-n/q_i$, that is,
		$\theta\neq0$. Fix $h\in\RR^n$ and $R>c|h|$. Since $\lambda R>c|\lambda h|$
		for every $\lambda>0$, the H\"ormander condition can be applied to the pair
		$(\lambda h,\lambda R)$, and \eqref{T-scaling} gives
		$\lambda^{\theta}\Lambda_\beta(h,R)\le C$ for all $\lambda>0$. Letting
		$\lambda\to\infty$ if $\theta>0$, or $\lambda\to0$ if $\theta<0$, we
		conclude that $\Lambda_\beta(h,R)=0$. Consequently $N_j(h,R)=0$ for every $j\ge1$,
		which means that $k_i(x-h)=k_i(x)$ for a.e. $x$ with $|x|>2R$.
		
		This holds for every $R>c|h|$, and we can let $R$ decrease to $c|h|$.
		Indeed, let $h\neq0$ (the case $h=0$ is trivial) and $R_l=c|h|(1+\frac1l)$.
		Since $\{|x|>2c|h|\}=\bigcup_l\{|x|>2R_l\}$ and a countable union of null
		sets is null, we obtain $k_i(x-h)=k_i(x)$ for a.e. $x$ with $|x|>2c|h|$.
		In other words, for every $r>0$ the function $k_i$ is invariant under all
		sufficiently small translations on the exterior region $U_r=\{x:|x|>r\}$:
		\begin{equation}\label{T-invariance}
			k_i(\cdot-h)=k_i\quad\text{a.e. in } U_r,
			\qquad\text{whenever } |h|<\frac r{2c}.
		\end{equation}
		
		We next show that this invariance forces $k_i$ to be locally constant on
		$U_r$. Since $\Omega_i\in L^1(S^{n-1})$, we have
		$k_i\in L^1_{loc}(\RR^n\setminus\{0\})$. Let $\varphi\in C_c^\infty(U_r)$
		and $1\le l\le n$. By \eqref{T-invariance}, for $|t|$ small enough,
		\[
		\int k_i(x)\frac{\varphi(x+te_l)-\varphi(x)}{t}\,dx
		=\int\frac{k_i(x-te_l)-k_i(x)}{t}\varphi(x)\,dx=0 ,
		\]
		and letting $t\to0$, by dominated convergence, we get
		$\int k_i\,\partial_l\varphi=0$. Thus $\nabla k_i=0$ in $\mathcal D'(U_r)$,
		and therefore $k_i$ coincides a.e. with a constant $c_V$ on each connected
		component $V$ of $U_r$.
		
		To identify these components, observe that the map
		$\Phi(x)=(x/|x|,|x|)$ is a homeomorphism from $U_r$ onto
		$S^{n-1}\times(r,\infty)$. Hence, if $n\ge2$, $U_r$ is connected, while if
		$n=1$ its connected components are $(r,\infty)$ and $(-\infty,-r)$. In both
		cases, since the dilation $x\mapsto2x$ corresponds under $\Phi$ to
		$(\theta,\rho)\mapsto(\theta,2\rho)$, every connected component $V$ of $U_r$
		satisfies $2V\subset V$. Using the homogeneity of $k_i$ once more, we
		obtain $c_V=k_i(2x)=2^{-n/q_i}k_i(x)=2^{-n/q_i}c_V$ for a.e. $x\in V$, and
		since $2^{-n/q_i}\ne1$, it follows that $c_V=0$. Hence $k_i=0$ a.e. on $U_r$
		for every $r>0$, that is, $k_i=0$ a.e. on $\RR^n$, which contradicts
		$k_i\not\equiv0$. This proves (a).
		
		Finally, if $\delta\ne0$, then $n-n/q_i-\delta\ne n-n/q_i$, so (b) is an
		immediate consequence of (a).
	\end{proof}
	
	\begin{remark}
		By Theorem~\ref{T-index}(c), $(\mathrm H_2)$ is compatible with the
		homogeneity of $k_i$: since $2^{j\delta}=(2^jR/R)^{\delta}$, the factor does
		not depend on the scale. Writing $2^{j\delta}=(2^jR)^{\delta}R^{-\delta}$,
		condition $(\mathrm H_2)$ reads
		\[
		\sum_{j=1}^{\infty}(2^jR)^{n/q_i+\delta}
		\|k_i(\cdot-h)-k_i\|_{p_i,|x|\sim2^jR}\le C_i\,R^{\delta},
		\qquad R>c_i|h| .
		\]
		The factor $R^{\delta}$ on the right-hand side is exactly what a classical
		H\"ormander condition does not allow and, by Theorem~\ref{T-index}(b), it
		cannot be removed by changing the index. On the other hand,
		Remark~\ref{R-counterexample} shows that the classical condition
		$H_{n-n/q_i,p_i}$ is not sufficient for our results, for any
		$1\le p_i\le\infty$.
	\end{remark}
	
	\section{Proof of the main result}\label{S-proof}
	
	In this section we prove Theorem~\ref{T-HpLq}. We first state an extension
	theorem and prove the technical lemma, which gives the uniform boundedness of
	$T_\alpha$ on atoms; then we obtain the unweighted estimate needed to apply
	them and prove the theorem. Finally, we show that $(\mathrm H_2)$ cannot be
	replaced by the classical H\"ormander condition.
	
	\subsection{An extension theorem}
	We state the extension result with two weights, which is the form needed
	for Theorem~\ref{T-HpLq}. It extends \cite[Theorem~2.10]{Rocha}, which
	corresponds to $q_0=p_0$, $u=v$ and $q=p$.
	
	\begin{theorem}\label{Extension}
		Let $T$ be a bounded linear operator from $L^{p_0}(\RR^{n})$ into
		$L^{q_0}(\RR^{n})$ for some $1<p_0\le q_0<\infty$. Let $v\in\mathcal A_\infty$,
		$0<p\le1$ with $p\frac{r_v}{r_v-1}<p_0$, $d\ge\lfloor n(\frac{\tilde q_v}{p}-1)\rfloor$,
		$u$ a weight, and $p\le q<\infty$. Then $T$ extends to a bounded operator
		from $H^{p}_v(\RR^{n})$ into $L^{q}_u(\RR^{n})$, that is, there is $C>0$
		such that $\|Tf\|_{L^q_u}\le C\|f\|_{H^p_v}$ for every $f\in\mathcal D_0$,
		if and only if
		\[
		\sup\bigl\{\|Ta\|_{L^q_u}: a \text{ is a } v\text{-}(p,p_0,d)\text{ atom}\bigr\}<\infty .
		\]
	\end{theorem}
	
	\begin{proof}
		Assume first that $\|Tf\|_{L^q_u}\le C_1\|f\|_{H^p_v}$ for every
		$f\in\mathcal D_0$, and let $a$ be a $v$-$(p,p_0,d)$ atom. Since $a$ need not
		belong to $\mathcal D_0$, we approximate it. Fix $\zeta\in\mathcal S(\RR^n)$
		such that $\widehat\zeta\in C_c^\infty(\RR^n)$ and $\widehat\zeta\equiv1$ in a
		neighborhood of the origin, let $\zeta_\lambda(x)=\lambda^{-n}\zeta(x/\lambda)$
		for $\lambda>0$, and set
		\[
		f_k=a*\zeta_{1/k}-a*\zeta_k,\qquad k\in\mathbb N.
		\]
		Since $a$ has compact support, $\widehat a\in C^\infty(\RR^n)$, and
		$\widehat{f_k}(\xi)=\widehat a(\xi)\bigl(\widehat\zeta(\xi/k)-\widehat\zeta(k\xi)\bigr)$
		is smooth, has compact support and vanishes in a neighborhood of the origin;
		hence $f_k\in\mathcal D_0$. As $k\to\infty$, $a*\zeta_{1/k}\to a$ in
		$L^{p_0}$, because $\int\zeta=\widehat\zeta(0)=1$, while
		\[
		\|a*\zeta_k\|_{p_0}\le\|a\|_1\|\zeta_k\|_{p_0}
		=k^{-n(1-1/p_0)}\|a\|_1\|\zeta\|_{p_0}\longrightarrow0 ;
		\]
		thus $f_k\to a$ in $L^{p_0}$.
		
		The norms $\|f_k\|_{H^p_v}$ are uniformly bounded. Indeed, let $\psi\in S_N$,
		$t,\lambda>0$ and $\sigma=\max\{t,\lambda\}$. Then
		$\psi_t*\zeta_\lambda=(\psi_{t/\sigma}*\zeta_{\lambda/\sigma})_\sigma$, and
		either $t/\sigma=1$ or $\lambda/\sigma=1$. A direct computation shows that the
		functions $\psi*\zeta_s$ and $\psi_s*\zeta$, with $0<s\le1$ and
		$\psi\in S_N$, are uniformly bounded in the semi-norms
		$\|\cdot\|_{\mu,\kappa}$ with $|\mu|,|\kappa|\le N-n-1$. Hence
		$M_N(g*\zeta_\lambda)\le C\,M_{N-n-1}g$ with $C$ independent of $\lambda$
		and, since the space $H^p_v$ does not depend on $N$ for $N$ large
		(see \cite[Chapter~VI]{ST}), $\|g*\zeta_\lambda\|_{H^p_v}\le C\|g\|_{H^p_v}$ for every
		$g\in H^p_v(\RR^n)$ and $\lambda>0$. By Lemma~\ref{Lemma8},
		$\sup_k\|f_k\|_{H^p_v}\le C\|a\|_{H^p_v}\le C$.
		
		Finally, since $T$ is bounded from $L^{p_0}$ into $L^{q_0}$, $Tf_k\to Ta$ in
		$L^{q_0}$, so a subsequence converges a.e., and Fatou's lemma gives
		\[
		\|Ta\|_{L^q_u}\le\liminf_{k\to\infty}\|Tf_k\|_{L^q_u}
		\le C_1\sup_k\|f_k\|_{H^p_v}\le C,
		\]
		with $C$ independent of $a$.
		
		Conversely, assume $\|Ta\|_{L^q_u}\le C_0$ for all $v$-$(p,p_0,d)$ atoms.
		Let $f\in\mathcal D_0$ and let $f=\sum_j\lambda_ja_j$ be the decomposition
		of Theorem~\ref{Theorem9}, which converges in $L^{p_0}$. Since $T$ is
		bounded from $L^{p_0}$ into $L^{q_0}$, $\sum_j\lambda_jTa_j$ converges to
		$Tf$ in $L^{q_0}$, so a subsequence of the partial sums converges a.e. and
		\begin{equation}\label{Tpuntual}
			|Tf(x)|\le\sum_j|\lambda_j|\,|Ta_j(x)|\qquad\text{a.e. } x\in\RR^{n}.
		\end{equation}
		If $p\le q<1$, then \eqref{Tpuntual} and the inequality
		$(\sum_j b_j)^q\le\sum_j b_j^q$ give
		\[
		\|Tf\|_{L^q_u}^q\le\sum_j|\lambda_j|^q\|Ta_j\|^q_{L^q_u}
		\le C_0^q\sum_j|\lambda_j|^q,
		\]
		and $(\sum_j|\lambda_j|^q)^{1/q}\le(\sum_j|\lambda_j|^p)^{1/p}$ since
		$p\le q$. If $q\ge1$, Minkowski's inequality and monotone convergence give
		\[
		\|Tf\|_{L^q_u}\le\sum_j|\lambda_j|\,\|Ta_j\|_{L^q_u}
		\le C_0\sum_j|\lambda_j|\le C_0\Bigl(\sum_j|\lambda_j|^p\Bigr)^{1/p},
		\]
		since $p\le1$. In both cases $\|Tf\|_{L^q_u}\le C\|f\|_{H^p_v}$ for all
		$f\in\mathcal D_0$, and the result follows from Theorem~\ref{T-dense}, since
		$v$ is doubling.
	\end{proof}
	
	\subsection{The technical lemma}
	We now prove that $T_\alpha$ maps $w^{p}$-$(p,p_0,0)$ atoms uniformly into
	$L^q_{w^q}$. For a single kernel, an argument of this type is given in the
	proof of \cite[Theorem~5.10]{HW}, following \cite{DLeLi}; the main new
	difficulty here is the presence of several kernels composed with the
	matrices $A_i$, whose singularities are located at the different points
	$A_ix_0$.
	
	\begin{lemma}\label{TL}
		Let $A_1,\dots,A_m$ be invertible matrices, and let $T_\alpha$ satisfy
		$(\mathrm H_1)$ and $(\mathrm H_2)$,
		where $0<\alpha<n$ and $p,\delta,q$ are as in Theorem~\ref{T-HpLq}. Let
		$s_0=a(t,q)$ be as in \eqref{def-a}, let $\eta>q$, and assume
		\[
		w^{\max\{s_0,\eta\}}\in\mathcal A_1\qquad\text{and}\qquad
		w(A_ix)\le Cw(x)\ \text{a.e.},\ 1\le i\le m .
		\]
		Let $1<p_0<\frac n\alpha$, $\frac1{q_0}=\frac1{p_0}-\frac\alpha n$, be such
		that $T_\alpha:L^{p_0}\to L^{q_0}$ is bounded and
		\[
		q_0>q,\qquad q\Bigl(\frac{q_0}{q}\Bigr)'\le\eta .
		\]
		Then there is $C>0$ such that $\|T_\alpha a\|_{L^q_{w^q}}\le C$ for every
		function $a$ supported in a ball $B$ with
		$\|a\|_{p_0}\le|B|^{1/p_0}[w^p(B)]^{-1/p}$ and $\int a=0$; in particular,
		for every $w^p$-$(p,p_0,d)$ atom.
	\end{lemma}
	
	\begin{proof}
		Let $B=B(x_0,r)$ be the ball associated with $a$. The idea of the proof is
		to split $\RR^n$ into a region near the images $A_i(B)$ of the support of
		$a$, where we only use the size of $a$ and the $L^{p_0}-L^{q_0}$
		boundedness of $T_\alpha$, and a far region, where we use the cancellation
		of $a$ together with the regularity of the kernels given by
		$(\mathrm H_2)$. Before doing so, we collect the estimates on the weight and
		on the kernels that will be needed.
		
		We first note that, by H\"older's inequality and the size condition on $a$,
		\begin{equation}\label{TL-L1}
			\|a\|_1\le|B|^{1-1/p_0}\|a\|_{p_0}\le|B|\,[w^p(B)]^{-1/p}.
		\end{equation}
		Let
		\[
		\rho=\Bigl(\frac{q_0}{q}\Bigr)',\qquad
		\gamma_i=\frac n{q_i}-\frac n{p_i}>0,\quad 1\le i\le m .
		\]
		Since $p<q\le s_0$ and $q<q\rho\le\eta$, all the exponents $s_0$, $q\rho$,
		$p$ and $q$ are at most $\max\{s_0,\eta\}$, so Lemma~\ref{L-A1}(i), applied to
		$v=w^{\max\{s_0,\eta\}}$, shows that $w^{s_0},w^{q\rho},w^p,w^q\in\mathcal A_1$.
		In particular, Lemma~\ref{L-A1}(iii), applied to $v=w^{q\rho}$ with
		$\theta=1/\rho$, gives $w^q\in RH_\rho$, and, by Lemma~\ref{L-A1}(iv), the
		measures
		$w^p(x)\,dx$ and $w^q(x)\,dx$ satisfy $w^p(\lambda Q)\le C\lambda^n w^p(Q)$ and
		$w^q(\lambda Q)\le C\lambda^n w^q(Q)$ for every ball $Q$ and $\lambda\ge1$.
		
		The weights $w^p$ and $w^q$ can be compared on balls. Indeed, since
		$w^q\in\mathcal A_1$ and, trivially, $\einf_Q w\le(w^p(Q)/|Q|)^{1/p}$, we
		have
		\[
		\Bigl(\frac{w^q(Q)}{|Q|}\Bigr)^{1/q}\le C\einf_Q w
		\le C\Bigl(\frac{w^p(Q)}{|Q|}\Bigr)^{1/p},
		\]
		and, as $\frac1q-\frac1p=-\frac\alpha n$, this gives
		\begin{equation}\label{TL-comparison}
			[w^q(Q)]^{1/q}\le C|Q|^{-\alpha/n}[w^p(Q)]^{1/p}
			\qquad\text{for every ball } Q.
		\end{equation}
		In the same way, the condition $w^{s_0}\in\mathcal A_1$ allows us to
		control the $L^{s_0}$ norm of $w$ on subsets of a ball by its $w^q$ measure:
		for every measurable $E\subset Q$,
		\begin{equation}\label{TL-weight}
			\|w\chi_E\|_{s_0}\le\bigl(C|Q|\einf_Q w^{s_0}\bigr)^{1/s_0}
			\le C|Q|^{\frac1{s_0}-\frac1q}[w^q(Q)]^{1/q}.
		\end{equation}
		As for the kernels, $(\mathrm H_1)$ and Lemma~\ref{L-tail} give the tail
		estimate
		\begin{equation}\label{TL-tail}
			\|k_i\chi_{\{|x|>u\}}\|_{p_i}\le Cu^{-\gamma_i},\qquad u>0,\ 1\le i\le m.
		\end{equation}
		
		We now define the two regions. Let $c_i$ be the constants in
		$(\mathrm H_2)$, and fix
		\[
		M>\max\Bigl\{1,\max_{1\le i\le m}\|A_i\|,\max_{1\le i\le m}c_i\|A_i\|\Bigr\},
		\qquad
		L>\max\Bigl\{1,4M\max_{1\le i\le m}\|A_i^{-1}\|\Bigr\}.
		\]
		Let $B^*=B(x_0,Lr)$, $U=\bigcup_{i=1}^{m}A_i(B^*)$ and $D=\RR^n\setminus U$.
		Then
		\[
		\|T_\alpha a\|_{L^q_{w^q}}\le\sum_{i=1}^m I_i+I_{m+1},
		\qquad
		I_i=\|T_\alpha a\,\chi_{A_i(B^*)}\|_{L^q_{w^q}},\quad
		I_{m+1}=\|T_\alpha a\,\chi_D\|_{L^q_{w^q}} ,
		\]
		and it suffices to show that each of these terms is bounded by a constant
		independent of $a$.
		
		We begin with the terms $I_i$, $1\le i\le m$, which correspond to the near
		region. H\"older's inequality with exponents $q_0/q$ and $\rho$ gives
		\[
		I_i^q\le\|T_\alpha a\|_{q_0}^q
		\Bigl(\int_{A_i(B^*)}w^{q\rho}\Bigr)^{1/\rho}.
		\]
		To estimate the integral of the weight, we change variables $y=A_ix$ and
		use the hypothesis $w(A_ix)\le Cw(x)$, which reduces the integral over
		$A_i(B^*)$ to an integral over $B^*$; then the reverse H\"older property
		$w^q\in RH_\rho$ yields
		\begin{align*}
			\Bigl(\int_{A_i(B^*)}w^{q\rho}\Bigr)^{1/\rho}
			&=|\det A_i|^{1/\rho}\Bigl(\int_{B^*}w(A_ix)^{q\rho}dx\Bigr)^{1/\rho}\\
			&\le C\Bigl(\int_{B^*}w^{q\rho}\Bigr)^{1/\rho}
			\le C|B^*|^{\frac1\rho-1}w^q(B^*).
		\end{align*}
		On the other hand, $\|T_\alpha a\|_{q_0}\le C\|a\|_{p_0}$, and the size
		condition on $a$ together with the doubling property
		$w^p(B^*)\le Cw^p(B)$ and $|B|\le|B^*|$ lead to
		\[
		I_i^q\le C|B|^{\frac q{p_0}}[w^p(B)]^{-\frac qp}|B^*|^{\frac1\rho-1}w^q(B^*)
		\le C|B^*|^{\frac q{p_0}+\frac1\rho-1}[w^p(B^*)]^{-\frac qp}w^q(B^*).
		\]
		Since $\frac1\rho=1-\frac q{q_0}$, the exponent of $|B^*|$ is
		$\frac q{p_0}+\frac1\rho-1=q(\frac1{p_0}-\frac1{q_0})=\frac{q\alpha}n$, which
		is exactly compensated by \eqref{TL-comparison} applied with $Q=B^*$:
		\[
		I_i^q\le C|B^*|^{\frac{q\alpha}{n}}[w^p(B^*)]^{-\frac{q}{p}}
		|B^*|^{-\frac{q\alpha}{n}}[w^p(B^*)]^{\frac{q}{p}}=C .
		\]
		
		We now turn to $I_{m+1}$. Here $y$ is far from all the points $A_ix_0$,
		and we organize $D$ according to the nearest of these points and to the
		dyadic distance to it. If $y\in D$, then $|A_\ell^{-1}y-x_0|\ge Lr$, hence
		\[
		|y-A_\ell x_0|\ge\frac{Lr}{\|A_\ell^{-1}\|}>4Mr
		\qquad\text{for every } 1\le\ell\le m .
		\]
		For $1\le\ell\le m$ let $D_\ell$ be the set of points $y\in D$ such that
		$|y-A_\ell x_0|<|y-A_ix_0|$ for $i<\ell$ and
		$|y-A_\ell x_0|\le|y-A_ix_0|$ for $i>\ell$; that is, each point of $D$ is
		assigned to a nearest center $A_\ell x_0$, ties being assigned to the
		smallest index. The sets $D_\ell$ are measurable and form a partition of
		$D$, and so do the sets
		\[
		D_{\ell,j}=D_\ell\cap\{y:2^jMr<|y-A_\ell x_0|\le2^{j+1}Mr\},\qquad
		1\le\ell\le m,\ j\ge1 .
		\]
		If $y\in D_{\ell,j}$, $1\le i\le m$ and $z\in B$, then
		$|y-A_ix_0|\ge|y-A_\ell x_0|>2^jMr$, since $A_\ell x_0$ is a nearest center,
		and therefore, as $M>\|A_i\|$,
		\begin{equation}\label{TL-geometry}
			|y-A_iz|\ge|y-A_ix_0|-\|A_i\|\,|z-x_0|>(2^j-1)Mr\ge2^{j-1}Mr .
		\end{equation}
		In other words, on $D_{\ell,j}$ all the singularities of the kernel are at
		distance comparable to $2^jr$ or larger.
		
		Since $\int a=0$, we can write
		\[
		T_\alpha a(y)=\int_B[K(y,z)-K(y,x_0)]a(z)\,dz,
		\]
		and, writing the difference of the products as a telescoping sum,
		$|K(y,z)-K(y,x_0)|\le\sum_{i_0=1}^m F_{i_0,z}(y)$, where
		\[
		F_{i_0,z}(y)=\prod_{i<i_0}|k_i(y-A_ix_0)|\,
		|k_{i_0}(y-A_{i_0} z)-k_{i_0}(y-A_{i_0} x_0)|\prod_{i>i_0}|k_i(y-A_iz)| .
		\]
		In each $F_{i_0,z}$ only the $i_0$-th factor is a difference, and it will
		be controlled by $(\mathrm H_2)$, while the remaining factors will be
		controlled by the tail estimate \eqref{TL-tail}. To make this precise, let
		$h_{i_0}(z)=A_{i_0}(z-x_0)$ and
		\[
		N_{i_0,k}(z)=\bigl\|[k_{i_0}(\cdot-h_{i_0}(z))-k_{i_0}]
		\chi_{\{2^kMr<|\cdot|\le2^{k+1}Mr\}}\bigr\|_{p_{i_0}},\qquad k\ge1.
		\]
		Since $|h_{i_0}(z)|<\|A_{i_0}\|r$ and $M>c_{i_0}\|A_{i_0}\|$, we have
		$Mr>c_{i_0}|h_{i_0}(z)|$, so $(\mathrm H_2)$ can be applied with $R=Mr$ and
		$h=h_{i_0}(z)$. Taking into account that
		$\|f\chi_{\{|x|\sim s\}}\|_{p_{i_0}}=|B(0,2s)|^{1/p_{i_0}}\|f\|_{p_{i_0},|x|\sim s}$,
		it gives
		\begin{equation}\label{TL-Hsum}
			\sum_{k=1}^\infty2^{k\delta}(2^kMr)^{\gamma_{i_0}}N_{i_0,k}(z)\le C,
			\qquad\text{uniformly in } z\in B.
		\end{equation}
		
		We now estimate $F_{i_0,z}$ on each set $D_{\ell,j}$. By \eqref{TL-geometry}
		and \eqref{TL-tail}, for $i\neq i_0$ the $L^{p_i}$ norm on $D_{\ell,j}$ of
		the factor $k_i(\cdot-A_ix_0)$ or $k_i(\cdot-A_iz)$ is bounded by
		$C(2^jr)^{-\gamma_i}$. For the difference factor, the change of variables
		$x=y-A_{i_0} x_0$ and the inclusion $D_{\ell,j}-A_{i_0} x_0\subset\{|x|>2^jMr\}$
		show that
		\[
		\bigl\|[k_{i_0}(\cdot-A_{i_0} z)-k_{i_0}(\cdot-A_{i_0} x_0)]\chi_{D_{\ell,j}}\bigr\|_{p_{i_0}}
		\le\sum_{k=j}^\infty N_{i_0,k}(z).
		\]
		For the weight, let $Q_\ell=B(A_\ell x_0,2Mr)$. Since
		$D_{\ell,j}\subset 2^j\overline{Q_\ell}$, \eqref{TL-weight} and the doubling
		property of $w^q$ yield
		\[
		\|w\chi_{D_{\ell,j}}\|_{s_0}
		\le C(2^jr)^{\frac n{s_0}-\frac nq}[w^q(2^jQ_\ell)]^{1/q}
		\le C(2^jr)^{\frac n{s_0}-\frac nq}2^{\frac{jn}{q}}[w^q(Q_\ell)]^{1/q}.
		\]
		We combine these three estimates by means of the generalized H\"older
		inequality with $\frac1q=\sum_{i=1}^m\frac1{p_i}+\frac1{s_0}$. Since, by
		\eqref{q_i} and \eqref{def-a},
		\[
		-\sum_{i\ne i_0}\gamma_i+\frac n{s_0}-\frac nq=\alpha-n+\gamma_{i_0}
		\qquad\text{and}\qquad \alpha-n+\frac nq=\delta ,
		\]
		the powers of $2^j r$ combine into
		\begin{equation}\label{TL-annular}
			\|F_{i_0,z}\chi_{D_{\ell,j}}\|_{L^q_{w^q}}
			\le Cr^{\alpha-n}[w^q(Q_\ell)]^{1/q}\,2^{j\delta}(2^jr)^{\gamma_{i_0}}
			\sum_{k=j}^\infty N_{i_0,k}(z).
		\end{equation}
		It is at this point that the factor $2^{j\delta}$ appears, and it is
		absorbed by $(\mathrm H_2)$ when we sum in $j$. Indeed, since
		$\gamma_{i_0}+\delta>0$, interchanging the order of summation and using
		\eqref{TL-Hsum} we obtain
		\begin{align*}
			\sum_{j=1}^\infty2^{j\delta}(2^jr)^{\gamma_{i_0}}\sum_{k=j}^\infty N_{i_0,k}(z)
			&=r^{\gamma_{i_0}}\sum_{k=1}^\infty N_{i_0,k}(z)\sum_{j=1}^k2^{j(\gamma_{i_0}+\delta)}\\
			&\le C\sum_{k=1}^\infty2^{k\delta}(2^kr)^{\gamma_{i_0}}N_{i_0,k}(z)\le C,
		\end{align*}
		where the powers of $M$ have been absorbed into $C$. Therefore, since
		$q>1$, Minkowski's integral inequality and \eqref{TL-annular} give
		\[
		I_{m+1}\le\sum_{\ell=1}^m\sum_{j=1}^\infty\sum_{i_0=1}^m
		\int_B|a(z)|\,\|F_{i_0,z}\chi_{D_{\ell,j}}\|_{L^q_{w^q}}\,dz
		\le Cr^{\alpha-n}\|a\|_1\sum_{\ell=1}^m[w^q(Q_\ell)]^{1/q}.
		\]
		It remains to compare $w^q(Q_\ell)$ with $w^q(B)$. Since
		$Q_\ell\subset A_\ell\bigl(B(x_0,2M\|A_\ell^{-1}\|r)\bigr)$, the change of
		variables $y=A_\ell x$, the hypothesis $w(A_\ell x)\le Cw(x)$, and the
		doubling property of $w^q$ (note that
		$2M\|A_\ell^{-1}\|>2\|A_\ell\|\,\|A_\ell^{-1}\|\ge2$) give
		$w^q(Q_\ell)\le Cw^q(B)$. Finally, by \eqref{TL-L1} and
		\eqref{TL-comparison} with $Q=B$,
		\[
		I_{m+1}\le Cr^{\alpha-n}|B|[w^p(B)]^{-1/p}[w^q(B)]^{1/q}
		\le Cr^{\alpha-n}|B|^{1-\alpha/n}\le C .
		\]
		This, together with the bound for $I_1,\dots,I_m$, proves the lemma.
	\end{proof}
	
	\begin{remark}\label{R-origin}
		In the proof of Lemma~\ref{TL}, the factor $2^{j\delta}$ in
		\eqref{TL-annular} comes exactly from the doubling estimate
		$w^q(2^jQ_\ell)\le C2^{jn}w^q(Q_\ell)$, combined with the factor
		$r^{\alpha-n}$. Indeed, after the generalized H\"older inequality, the
		powers of $2^j$ that are not absorbed into $(2^jr)^{\gamma_{i_0}}$ are
		$2^{j(\alpha-n)}$ and $2^{jn/q}$, and $\alpha-n+n/q=\delta$. For
		$w\equiv1$ this doubling estimate is an identity,
		$|2^jQ_\ell|=2^{jn}|Q_\ell|$, so no loss occurs at this step. This explains
		why the counterexample in Remark~\ref{R-counterexample} already works
		without weights.
	\end{remark}
	
	\subsection{The unweighted estimate}
	The following lemma is a consequence of the results in \cite{FR}.
	
	\begin{lemma}\label{L-unweighted}
		Assume $(\mathrm H)$ and $(\mathrm H_1)$, and that
		$k_i\in H_{n-n/q_i,p_i}$ for $1\le i\le m$ (this holds, for instance, under
		$(\mathrm H_2)$, by Remark~\ref{R-inclusions}). Let $t$ be as in
		\eqref{def-a}. If $t<p_0<\frac n\alpha$ and
		$\frac1{q_0}=\frac1{p_0}-\frac\alpha n$, then $T_\alpha$ is bounded from
		$L^{p_0}(\RR^n)$ into $L^{q_0}(\RR^n)$.
	\end{lemma}
	
	\begin{proof}
		Set $\alpha_i=n/q_i$, so that $0<\alpha_i<n$ and
		$\alpha_1+\dots+\alpha_m=n-\alpha$ by \eqref{q_i}. We first reduce the
		problem to finite exponents, since the results of \cite{FR} are stated for
		Young functions. Let $r_i=p_i$ if $p_i<\infty$, and let $r_i<\infty$ be large
		if $p_i=\infty$. By Remark~\ref{R-inclusions},
		$k_i\in S_{n-n/q_i,r_i}\cap H_{n-n/q_i,r_i}$ for $1\le i\le m$. Let
		$\tilde t$ be defined by $\frac1{r_1}+\dots+\frac1{r_m}+\frac1{\tilde t}=1$.
		Then $t\le\tilde t$, and $\tilde t\to t$ as the exponents $r_i$
		corresponding to $p_i=\infty$ tend to infinity; since $t<p_0$, we may choose
		them so that $\tilde t<p_0$.
		
		Let $\Psi_i(u)=u^{r_i}$ and $\varphi(u)=u^{\tilde t}$, which are Young
		functions. The normalized annular norms of \cite{FR} coincide with ours
		(the average is taken over $B(0,2s)$), so the conditions
		$k_i\in S_{n-n/q_i,r_i}\cap H_{n-n/q_i,r_i}$ mean that
		$k_i\in S_{n-\alpha_i,\Psi_i}\cap H_{n-\alpha_i,\Psi_i}$ in the notation of
		\cite[Section~2]{FR}. Since the constant $c$ in the H\"ormander condition
		may be enlarged, we may assume $c>1$, as required in \cite{FR}. Moreover,
		$\Psi_1^{-1}(u)\cdots\Psi_m^{-1}(u)\varphi^{-1}(u)=u$, so $\varphi$ is the
		complementary function of $\Psi_1,\dots,\Psi_m$, and
		$M_{\alpha,\varphi}=M_{\alpha,\tilde t}$.
		
		Let $f\in L^\infty_c(\RR^n)$. We first check that $T_\alpha f\in L^{q_0}$.
		If $m\ge2$, since $p_0>1$ we have $q_0>\frac n{n-\alpha}>\frac{n-\alpha}n$,
		so \cite[Lemma 5.1]{FR} with $\nu\equiv1$ gives $T_\alpha f\in L^{q_0}$. (Its
		proof shows that $|T_\alpha f(x)|\le C|x|^{\alpha-n}$ for $|x|$ large, and
		this function is in $L^{q_0}$ near infinity precisely because
		$q_0>\frac n{n-\alpha}$.) If $m=1$, then $T_\alpha f$ is bounded, since $k_1$
		is integrable near the origin by the size condition. Moreover, if
		$\operatorname{supp}f\subset B(0,\rho)$ and $|x|>2\|A_1\|\rho$, H\"older's
		inequality and $k_1\in S_{n-n/q_1,r_1}$ give
		$|T_\alpha f(x)|\le C|x|^{\alpha-n+n/r_1}$, and this function is in
		$L^{q_0}$ near infinity because
		$\frac1{q_0}=\frac1{p_0}-\frac\alpha n<\frac1{\tilde t}-\frac\alpha n
		=1-\frac\alpha n-\frac1{r_1}$, since $p_0>\tilde t=r_1'$. Hence
		\cite[Theorem 3.3]{FR} with $w\equiv1$ and exponent $q_0$ yields
		\[
		\|T_\alpha f\|_{q_0}\le C\|M_{\alpha,\tilde t}f\|_{q_0}
		= C\bigl\|M_{\alpha\tilde t}(|f|^{\tilde t})\bigr\|_{q_0/\tilde t}^{1/\tilde t}
		\le C\||f|^{\tilde t}\|_{p_0/\tilde t}^{1/\tilde t}=C\|f\|_{p_0},
		\]
		where we used \cite[Remark 2.2]{FR} and the boundedness of the fractional
		maximal operator $M_{\alpha\tilde t}$ from $L^{p_0/\tilde t}$ into
		$L^{q_0/\tilde t}$, which holds because
		$1<\frac{p_0}{\tilde t}<\frac n{\alpha\tilde t}$ and
		$\frac{\tilde t}{q_0}=\frac{\tilde t}{p_0}-\frac{\alpha\tilde t}n$. The
		result follows by density of $L^\infty_c$ in $L^{p_0}$.
	\end{proof}
	
	\begin{remark}
		Hypothesis $(\mathrm H)$ is not used in the proof of Lemma~\ref{TL}, which
		holds for arbitrary invertible matrices $A_1,\dots,A_m$; in the proof of
		Theorem~\ref{T-HpLq} it enters only through Lemma~\ref{L-unweighted}.
	\end{remark}
	
	\subsection{Proof of Theorem~\ref{T-HpLq}}
	Let $s_0=a(t,q)\ge q$. By Lemma~\ref{L-A1}(ii) applied to $w^{s_0}$, there
	is $\varepsilon>0$ such that $w^{\eta}\in\mathcal A_1$ with
	$\eta=s_0(1+\varepsilon)>q$. In particular $w^{\max\{s_0,\eta\}}=w^\eta\in\mathcal A_1$.
	
	Since $\eta>q>p$ and $t<\frac n\alpha$, we can choose $p_0$ with
	\[
	\max\Bigl\{1,\,t,\,\frac{p\eta}{\eta-p},\,
	\Bigl(\frac\alpha n+\frac1q-\frac1\eta\Bigr)^{-1}\Bigr\}
	<p_0<\frac n\alpha .
	\]
	The interval is nonempty: $\frac1p-\frac1\eta>\frac1p-\frac1q=\frac\alpha n$
	and $\frac\alpha n+\frac1q-\frac1\eta>\frac\alpha n$. With
	$\frac1{q_0}=\frac1{p_0}-\frac\alpha n$, the last lower bound gives
	$\frac1q-\frac1{q_0}>\frac1\eta$, that is, $q_0>q$ and
	$q(\frac{q_0}q)'<\eta$. Since $t<p_0<\frac n\alpha$, Lemma~\ref{L-unweighted}
	gives the boundedness of $T_\alpha:L^{p_0}\to L^{q_0}$.
	
	By Lemma~\ref{L-A1}(iii), applied to $v=w^\eta$ with $\theta=p/\eta$,
	$w^p\in RH_{\eta/p}$, so $r_{w^p}\ge\eta/p$ and
	\[
	p\frac{r_{w^p}}{r_{w^p}-1}\le\frac{p\eta}{\eta-p}<p_0 .
	\]
	Since $w^p\in\mathcal A_1$, we have $\tilde q_{w^p}=1$, and $\delta<1$
	gives $\lfloor n(\frac1p-1)\rfloor=0$; hence we may take $d=0$. By
	Lemma~\ref{TL}, $\|T_\alpha a\|_{L^q_{w^q}}\le C$ for every
	$w^p$-$(p,p_0,0)$ atom $a$. Since $q>1\ge p$, Theorem~\ref{Extension} with
	$v=w^p$ and $u=w^q$ shows that $T_\alpha$ extends to a bounded operator from
	$H^p_{w^p}(\RR^n)$ into $L^q_{w^q}(\RR^n)$. \qed
	
	\begin{remark}\label{R-HW}
		For a single kernel, the result closest to Theorem~\ref{T-HpLq} is
		\cite[Theorem~5.10]{HW}, which states the following: if $0<\alpha\le1$,
		$\frac n{n+\alpha}\le p<1$, $\frac1q=\frac1p-\frac\alpha n$,
		$s>\frac n{n-\alpha}$, the kernel $K_\alpha$ satisfies the $L^1$ size
		condition $\int_{r<|x|\le2r}|K_\alpha(x)|\,dx\le Cr^{\alpha}$ for $r>0$ and
		the fractional $L^s$-H\"ormander condition of order $\alpha$, that is,
		$K_\alpha\in H^{(\alpha)}_{\alpha,s}$, and
		$\omega^{ns/(ns-n-s\alpha)}\in\mathcal A_1$, then
		$f\mapsto\int K_\alpha(x-y)f(y)\,dy$ is bounded from
		$H^p_{\omega^p}(\RR^n)$ into $L^q_{\omega^q}(\RR^n)$.
		
		Theorem~\ref{T-HpLq} with $m=1$, $A_1=I$, $k_1=K_\alpha$,
		$q_1=\frac n{n-\alpha}$ and $p_1=s$ applies to the same operators. In this
		case, $(\mathrm H_1)$ and $(\mathrm H_2)$ mean that
		$K_\alpha\in S_{\alpha,s}\cap H^{(\delta)}_{\alpha,s}$, and \eqref{def-a}
		gives $a(t,q)=\frac{sq}{s-q}$. The main differences between both results are
		the following.
		\begin{itemize}
			\item \emph{Operators.} Theorem~\ref{T-HpLq} covers $m\ge2$ kernels
			composed with matrices $A_i$, which are not considered in \cite{HW}. Even for
			$m=1$, it allows a matrix $A_1\neq I$ and the case $s=\infty$.
			\item \emph{Smoothness.} Theorem~\ref{T-HpLq} requires the
			H\"ormander condition of order $\delta=n(\frac1p-1)$ instead of $\alpha$.
			Since $\delta<\alpha$ when $p>\frac n{n+\alpha}$, this is a weaker
			requirement and, for $n\ge2$, it is strictly weaker by
			Proposition~\ref{P-strict}.
			\item \emph{Range of $\alpha$.} In Theorem~\ref{T-HpLq}, $\alpha$ may
			be any number in $(0,n)$, since only $\delta<1$ is needed, while
			\cite[Theorem~5.10]{HW} requires $0<\alpha\le1$.
			\item \emph{Weights.} The condition on the weight in
			Theorem~\ref{T-HpLq} is weaker. Indeed, with $a=\frac{sq}{s-q}$,
			\[
			\frac1a=\frac1p-\frac\alpha n-\frac1s>1-\frac\alpha n-\frac1s
			=\frac{ns-n-s\alpha}{ns},
			\]
			so that $a<\frac{ns}{ns-n-s\alpha}$ and, by Lemma~\ref{L-A1}(i),
			$\omega^{ns/(ns-n-s\alpha)}\in\mathcal A_1$ implies
			$\omega^{a}\in\mathcal A_1$. This improvement comes from the
			self-improvement property of $\mathcal A_1$ weights (Lemma~\ref{L-A1}(ii))
			used in the proof of Theorem~\ref{T-HpLq}.
			\item \emph{Range of $p$ and size condition.} On the other hand,
			\cite[Theorem~5.10]{HW} includes the endpoint $p=\frac n{n+\alpha}$ and only
			requires an $L^1$ size condition, while Theorem~\ref{T-HpLq} requires
			$p>\frac n{n+\alpha}$ and the $L^s$ size condition $K_\alpha\in S_{\alpha,s}$.
		\end{itemize}
	\end{remark}
	
	\subsection{Sharpness of \texorpdfstring{$(\mathrm H_2)$}{(H2)}}\label{S-sharp}
	We now show that $(\mathrm H_2)$ cannot be replaced by the classical
	fractional H\"ormander condition.
	
	\begin{remark}[The factor $2^{j\delta}$ in $(\mathrm H_2)$ cannot be removed]
		\label{R-counterexample}
		Let $n\ge2$, $0<\alpha<n$, $p$ as in Theorem~\ref{T-HpLq}, and assume
		$\delta>1/q$. For instance, $n=2$, $\alpha=9/5$ and $p=10/13$ give $q=5/2$,
		$\delta=3/5$ and $\delta-1/q=1/5$. Fix $0<\varepsilon\le\delta-1/q$, let
		$m=2$, $q_1=q_2=2n/(n-\alpha)$, so that \eqref{q_i} holds, $A_1=I$,
		$A_2=-I$, $w\equiv1$, and
		\[
		k_1(x)=\frac{|x_n|^{\varepsilon}}{|x|^{n/q_1+\varepsilon}}
		=\frac{\Omega_1(x)}{|x|^{n/q_1}},\quad
		\Omega_1(x)=\Bigl|\frac{x_n}{|x|}\Bigr|^{\varepsilon},
		\qquad
		k_2(x)=\frac{1}{|x|^{n/q_2}} .
		\]
		
		\smallskip\noindent\textit{The hypotheses hold.}
		$(\mathrm H)$ holds, since $A_1$, $A_2$ and $A_1-A_2=2I$ are invertible.
		Since $\Omega_1,\Omega_2\in L^\infty(S^{n-1})$, Example~\ref{Ex-homog} gives
		$(\mathrm H_1)$ with $p_1=p_2=\infty$; then $t=1$ and $a(t,q)=q$. Clearly
		$w^q\equiv1\in\mathcal A_1$ and $w(A_ix)=w(x)$. Since
		$\bigl||s|^\varepsilon-|s'|^\varepsilon\bigr|\le|s-s'|^\varepsilon$,
		$\Omega_1$ is H\"older continuous of order $\varepsilon$ on $S^{n-1}$, and a
		standard computation gives, for $|x|\ge2|h|$,
		\[
		|k_1(x-h)-k_1(x)|\le C\Bigl(\frac{|h|}{|x|}\Bigr)^{\varepsilon}|x|^{-n/q_1},
		\qquad
		|k_2(x-h)-k_2(x)|\le C\frac{|h|}{|x|}|x|^{-n/q_2}.
		\]
		Consequently, for $0\le\delta'<\varepsilon$, $c_i=2$ and $R>2|h|$,
		\[
		\sum_{j=1}^\infty2^{j\delta'}(2^jR)^{n/q_i}
		\|k_i(\cdot-h)-k_i\|_{\infty,|x|\sim2^jR}
		\le C\Bigl(\frac{|h|}{R}\Bigr)^{\varepsilon}
		\sum_{j=1}^\infty2^{j(\delta'-\varepsilon)}\le C,
		\qquad i=1,2.
		\]
		In particular $k_i\in H_{n-n/q_i,\infty}$, and in fact
		$k_i\in H^{(\delta')}_{n-n/q_i,\infty}$ for every $\delta'<\varepsilon$.
		Moreover $T_\alpha$ is bounded from $L^{p_0}$ into $L^{q_0}$ for every
		$1<p_0<n/\alpha$, $\frac1{q_0}=\frac1{p_0}-\frac\alpha n$. This follows from
		Lemma~\ref{L-unweighted}, since $t=1$, and also from \cite{GU2}, since
		$|k_1(x)|\le|x|^{-n/q_1}$ gives
		$|T_\alpha f(x)|\le\int|x-y|^{-n/q_1}|x+y|^{-n/q_2}|f(y)|\,dy$.
		
		\smallskip\noindent\textit{An atom whose image is not in $L^q$.}
		Write $z=(\tilde z,z_n)\in\RR^{n-1}\times\RR$. Let
		$\varphi\in C_c^\infty(\RR^{n-1})$ with $\varphi\ge0$,
		$\operatorname{supp}\varphi\subset\{|\tilde z|<1/2\}$ and $\int\varphi=1$,
		let $\psi=\chi_{[0,1/2]}-\chi_{[-1/2,0)}$, and put
		$a(z)=\varphi(\tilde z)\psi(z_n)$. Then $a\in L^\infty$,
		$\operatorname{supp}a\subset B(0,1)$ and $\int a=0$. Since $w\equiv1$ and
		$\delta<1$, a suitable multiple of $a$ is a $(p,p_0,0)$-atom. Define
		\[
		G(u)=\int_{\RR}|u-s|^{\varepsilon}\psi(s)\,ds,\qquad u\in\RR.
		\]
		The function $G$ is continuous and
		\[
		G(1/2)=\int_0^{1/2}\tau^{\varepsilon}\,d\tau-\int_{1/2}^{1}\tau^{\varepsilon}\,d\tau<0,
		\]
		so there exist $c>0$ and an interval $I\subset(0,1)$ such that $|G|\ge c$
		on $I$. Let $y=(\tilde y,u)$ with $u\in I$ and $|\tilde y|=\rho$. For
		$z\in\operatorname{supp}a$, uniformly,
		\[
		|y-z|^{-n/q_1-\varepsilon}\,|y+z|^{-n/q_2}
		=\rho^{\alpha-n-\varepsilon}\bigl(1+O(\rho^{-1})\bigr),
		\]
		by \eqref{q_i}. Since $\int\varphi=1$,
		\[
		T_\alpha a(y)
		=\int|u-z_n|^{\varepsilon}|y-z|^{-n/q_1-\varepsilon}|y+z|^{-n/q_2}a(z)\,dz
		=\rho^{\alpha-n-\varepsilon}\bigl(G(u)+O(\rho^{-1})\bigr).
		\]
		Hence there is $\rho_1>0$ such that
		$|T_\alpha a(y)|\ge\frac c2\rho^{\alpha-n-\varepsilon}$ whenever
		$|\tilde y|=\rho\ge\rho_1$ and $u\in I$. Integrating in polar coordinates in
		$\tilde y$ (with counting measure on $S^0$ when $n=2$),
		\[
		\int_{\RR^n}|T_\alpha a|^q
		\ge\Bigl(\frac c2\Bigr)^q|I|\,|S^{n-2}|
		\int_{\rho_1}^\infty\rho^{(\alpha-n-\varepsilon)q+n-2}\,d\rho=\infty,
		\]
		because $(\alpha-n-\varepsilon)q+n-2\ge-1$ is equivalent to
		$\varepsilon\le\alpha-n+\frac{n-1}q=\delta-\frac1q$.
		
		\smallskip
		Therefore $T_\alpha a\notin L^q$. Hence the conclusion of Lemma~\ref{TL}
		fails and, by Theorem~\ref{Extension} with $v=u\equiv1$, $T_\alpha$ does not
		extend to a bounded operator from $H^p(\RR^n)$ into $L^q(\RR^n)$, so the
		conclusion of Theorem~\ref{T-HpLq} fails as well, although the kernels satisfy
		$(\mathrm H_1)$, the classical H\"ormander condition, and even
		$k_i\in H^{(\delta')}_{n-n/q_i,\infty}$ for every $\delta'<\varepsilon$.
		In particular, $\delta$ cannot be replaced in $(\mathrm H_2)$ by any
		$\delta'<\delta-1/q$.
		
		The same example gives more information when $(\mathrm H_1)$ and
		$(\mathrm H_2)$ are considered with a finite exponent for $k_1$. Let
		$q_1<p_1<\infty$ and $p_2=\infty$. Since $\Omega_1$ is bounded,
		$(\mathrm H_1)$ holds with these exponents. Moreover, since
		$0<\varepsilon\le\delta-\frac1q<1-\frac1{p_1}$, Proposition~\ref{P-strict}(a),
		with $\beta=n-n/q_1$ and $\tau=p_1$, gives
		$k_1\in H^{(\delta')}_{n-n/q_1,p_1}$ for every
		$\delta'<\varepsilon+\frac1{p_1}$, while
		$k_2\in H^{(\delta')}_{n-n/q_2,\infty}$ for every $\delta'<1$. Taking
		$\varepsilon=\delta-\frac1q$, we conclude that, for this choice of exponents,
		$\delta$ cannot be replaced in $(\mathrm H_2)$ by any
		$\delta'<\delta-\frac1q+\frac1{p_1}$. This is consistent with
		Theorem~\ref{T-HpLq}, since $\frac1{p_1}<\frac1{q_1}<1-\frac\alpha n<\frac1q$,
		so that $\delta-\frac1q+\frac1{p_1}<\delta$.
	\end{remark}
	
	\section{Examples of matrices and weights}\label{S-examples}
	
	In this section we give examples of matrices and weights satisfying the
	hypotheses of Theorem~\ref{T-HpLq}, that is, $(\mathrm H)$ for the matrices
	and
	\begin{equation}\label{hyp-w}
		w^{a(t,q)}\in\mathcal A_1,\qquad
		w(A_ix)\le C\,w(x)\ \text{ a.e. } x\in\RR^n,\ 1\le i\le m,
	\end{equation}
	for the weight. We also discuss the choice of the exponents $p_i$ and give a
	concrete family of operators.
	
	\subsection{Matrices}
	
	\begin{example}\label{Ex-matrices}
		The following families of matrices satisfy $(\mathrm H)$.
		\begin{itemize}
			\item[(a)] $A_i=a_iA$, where $A$ is an invertible matrix and
			$a_1,\dots,a_m$ are distinct nonzero real numbers, since
			$A_i-A_j=(a_i-a_j)A$. In particular, $A_i=a_iI$, which is the setting of
			\cite{GU2}.
			\item[(b)] Diagonal matrices $A_i=\mathrm{diag}(d_{i,1},\dots,d_{i,n})$ with
			$d_{i,l}\neq0$ and $d_{i,l}\neq d_{j,l}$ for $i\neq j$ and $1\le l\le n$.
			\item[(c)] In $\RR^2$, identified with $\mathbb C$, the matrices of the maps
			$z\mapsto c_iz$, where $c_1,\dots,c_m\in\mathbb C$ are distinct and nonzero.
			Indeed, $A_i-A_j$ is the matrix of $z\mapsto(c_i-c_j)z$, whose determinant
			is $|c_i-c_j|^2$. For instance, the rotations $A_i=R_{\theta_i}$ with angles
			$\theta_1,\dots,\theta_m$ distinct modulo $2\pi$.
		\end{itemize}
	\end{example}
	
	\begin{remark}
		In odd dimensions rotations behave differently: if $n$ is odd and
		$Q,Q'\in SO(n)$, then $Q-Q'$ is not invertible, so two rotations never
		satisfy $(\mathrm H)$. Indeed, $P=Q'^{-1}Q\in SO(n)$ and, since
		$P-I=P(I-P^{T})$,
		\[
		\det(P-I)=\det(I-P^{T})=\det(I-P)=(-1)^n\det(P-I)=-\det(P-I),
		\]
		so $\det(P-I)=0$ and $Q-Q'=Q'(P-I)$ is singular. On the other hand,
		$A_1=I$ and $A_2=-I$ satisfy $(\mathrm H)$ in every dimension.
	\end{remark}
	
	\subsection{Weights}
	We first observe that $\mathcal A_1$ is invariant under linear changes of
	variables.
	
	\begin{lemma}\label{L-A1linear}
		Let $v\in\mathcal A_1$ and let $A$ be an invertible matrix. Then
		$v\circ A\in\mathcal A_1$.
	\end{lemma}
	
	\begin{proof}
		Let $B=B(x_B,r)$ be a ball. The ellipsoid $A(B)$ is contained in the ball
		$B'=B(Ax_B,\|A\|r)$, and
		$|B'|/|A(B)|=\|A\|^n/|\det A|\le(\|A\|\,\|A^{-1}\|)^n$. Hence
		\begin{align*}
			\frac1{|B|}\int_B v(Ax)\,dx&=\frac1{|A(B)|}\int_{A(B)}v
			\le\frac{(\|A\|\,\|A^{-1}\|)^n}{|B'|}\int_{B'}v\\
			&\le C\einf_{B'}v\le C\einf_{A(B)}v=C\einf_{B}\,(v\circ A).
		\end{align*}
	\end{proof}
	
	\begin{proposition}\label{P-radial}
		Let $\varphi:(0,\infty)\to(0,\infty)$ be nonincreasing and such that
		$\varphi(r/2)\le C\varphi(r)$ for every $r>0$, and let $w(x)=\varphi(|x|)$.
		Then, for every invertible matrix $A$, there is $C_A>0$ such that
		$w(Ax)\le C_Aw(x)$ for all $x\neq0$.
	\end{proposition}
	
	\begin{proof}
		Since $|Ax|\ge|x|/\|A^{-1}\|$ and $\varphi$ is nonincreasing,
		$w(Ax)\le\varphi(|x|/\|A^{-1}\|)$. If $\|A^{-1}\|\le1$, this is at most
		$\varphi(|x|)=w(x)$. Otherwise, let $k\in\mathbb N$ with
		$2^k\ge\|A^{-1}\|$; then $w(Ax)\le\varphi(2^{-k}|x|)\le C^kw(x)$.
	\end{proof}
	
	\begin{example}\label{Ex-radialweights}
		Let $a=a(t,q)$ and $-n<\gamma\le0$. The weights
		\[
		w_1(x)=|x|^{\gamma/a},\qquad w_2(x)=(1+|x|)^{\gamma/a}
		\]
		satisfy \eqref{hyp-w} for arbitrary invertible matrices $A_1,\dots,A_m$.
		Indeed, $w_1^a=|x|^\gamma\in\mathcal A_1$, and $w_2^a=(1+|x|)^\gamma$ is
		comparable to $\min\{1,|x|^\gamma\}$, which belongs to $\mathcal A_1$
		because $M(\min\{u,v\})\le\min\{Mu,Mv\}$. Moreover, both weights are of the
		form $\varphi(|x|)$ with $\varphi$ nonincreasing and
		$\varphi(r/2)\le2^{-\gamma/a}\varphi(r)$, so the second condition in
		\eqref{hyp-w} follows from Proposition~\ref{P-radial}. If the matrices $A_i$
		are orthogonal, every radial weight $w$ with $w^a\in\mathcal A_1$ satisfies
		\eqref{hyp-w}, with $w(A_ix)=w(x)$.
	\end{example}
	
	Non-radial examples can be obtained by symmetrization when the matrices
	generate a finite group.
	
	\begin{proposition}\label{P-symmetrization}
		Let $G$ be a finite group of invertible matrices containing
		$A_1,\dots,A_m$, let $v\in\mathcal A_1$, and set
		\[
		w=\Bigl(\sum_{g\in G}v\circ g\Bigr)^{1/a(t,q)}.
		\]
		Then $w$ satisfies \eqref{hyp-w}, and in fact $w(A_ix)=w(x)$ for a.e. $x$
		and $1\le i\le m$.
	\end{proposition}
	
	\begin{proof}
		By Lemma~\ref{L-A1linear}, $v\circ g\in\mathcal A_1$ for every $g\in G$, and
		a finite sum of $\mathcal A_1$ weights is an $\mathcal A_1$ weight, since
		$M(u_1+u_2)\le Mu_1+Mu_2$. Hence $w^{a(t,q)}\in\mathcal A_1$; since
		$a(t,q)\ge1$, $w$ is locally integrable. Moreover, for each $i$ the map
		$g\mapsto gA_i$ is a bijection of $G$, so
		$\sum_{g\in G}v(gA_ix)=\sum_{g\in G}v(gx)$, that is, $w(A_ix)=w(x)$.
	\end{proof}
	
	\begin{example}\label{Ex-symm}
		Let $A_1=I$ and $A_2=-I$, so that $G=\{I,-I\}$, and let
		$x_0\in\RR^n\setminus\{0\}$ and $-n<\gamma<0$. By
		Proposition~\ref{P-symmetrization} with $v(x)=|x-x_0|^\gamma$, the weight
		\[
		w(x)=\bigl(|x-x_0|^{\gamma}+|x+x_0|^{\gamma}\bigr)^{1/a(t,q)},
		\]
		which has two singularities, satisfies \eqref{hyp-w}. Similarly, in
		$\RR^2\simeq\mathbb C$, if $A_iz=c_iz$ with $c_1,\dots,c_m$ distinct $k$-th
		roots of unity and $G$ is the group of the $k$-th roots of unity, one
		obtains weights with $k$ singularities at the vertices of a regular polygon.
	\end{example}
	
	\begin{remark}\label{R-centered}
		For power weights with a singularity away from the origin the
		symmetrization cannot be avoided. Indeed, let $m\ge2$, assume
		$(\mathrm H)$, and let $w(x)=|x-x_0|^{\gamma}$ with $\gamma<0$. If
		$w(A_ix)\le Cw(x)$ a.e., then $A_ix_0=x_0$: otherwise $y_0=A_i^{-1}x_0\neq
		x_0$, and for $x$ in a small ball centered at $y_0$ the quotient
		$w(A_ix)/w(x)$ is arbitrarily large. Thus $A_1x_0=A_2x_0=x_0$, so
		$(A_1-A_2)x_0=0$ and $x_0=0$ by $(\mathrm H)$.
	\end{remark}
	
	\subsection{The choice of the exponents \texorpdfstring{$p_i$}{p\_i}}
	
	\begin{remark}\label{R-choice-pi}
		If $(\mathrm H_1)$ and $(\mathrm H_2)$ hold for some exponents $p_i$, the
		condition on the weight depends on them through
		$\frac1{a(t,q)}=\frac1q-\sum_i\frac1{p_i}$. Increasing the $p_i$ decreases
		$a(t,q)$ and, by Lemma~\ref{L-A1}(i), weakens the condition
		$w^{a(t,q)}\in\mathcal A_1$; the weakest one, $w^q\in\mathcal A_1$,
		corresponds to $p_i=\infty$ for every $i$. On the other hand, by
		Remark~\ref{R-inclusions}, $(\mathrm H_1)$ and $(\mathrm H_2)$ become
		stronger as the $p_i$ increase. For instance, if $\Omega_i$ is H\"older
		continuous of order $\varepsilon>\delta$ on $S^{n-1}$, one may take
		$p_i=\infty$ (Proposition~\ref{P-Holder}), while for a rough
		$\Omega_i\in L^{\tau}(S^{n-1})\setminus L^\infty(S^{n-1})$ satisfying
		\eqref{Dini-delta} one takes $p_i=\tau$ (Proposition~\ref{P-Dini}), at the
		price of a stronger condition on the weight.
	\end{remark}
	
	\subsection{A concrete family}
	
	\begin{example}\label{Ex-family}
		Let $n=2$, $\alpha=9/5$, $m=2$, $q_1=q_2=20$, $A_1=I$, $A_2=-I$ and, for
		$0<\varepsilon\le1$,
		\[
		k_1(x)=\Bigl|\frac{x_2}{|x|}\Bigr|^{\varepsilon}|x|^{-1/10},\qquad
		k_2(x)=|x|^{-1/10},
		\]
		which is the setting of Section~\ref{S-sharp} with $n=2$. Let $p=10/13$, so
		that $q=5/2$ and $\delta=3/5$.
		\begin{itemize}
			\item[(i)] If $3/5<\varepsilon\le1$, then $T_\alpha$ is bounded from
			$H^{10/13}_{w^{10/13}}(\RR^2)$ into $L^{5/2}_{w^{5/2}}(\RR^2)$ for every
			weight $w$ such that $w^{5/2}\in\mathcal A_1$ and $w(-x)\le Cw(x)$; for
			instance, $w(x)=|x|^{\gamma}$ or $w(x)=(1+|x|)^{\gamma}$ with
			$-4/5<\gamma\le0$, or $w(x)=(|x-x_0|^{\gamma}+|x+x_0|^{\gamma})^{2/5}$ with
			$x_0\neq0$ and $-2<\gamma<0$.
			\item[(ii)] If $0<\varepsilon\le1/5$, this fails even for $w\equiv1$.
		\end{itemize}
		Indeed, $\Omega_1$ and $\Omega_2\equiv1$ are bounded, so $(\mathrm H_1)$
		holds with $p_1=p_2=\infty$; then $t=1$ and $a(t,q)=q=5/2$. As in
		Section~\ref{S-sharp}, for $|x|\ge2|h|$,
		\[
		|k_1(x-h)-k_1(x)|\le C\Bigl(\frac{|h|}{|x|}\Bigr)^{\varepsilon}|x|^{-1/10},
		\qquad
		|k_2(x-h)-k_2(x)|\le C\frac{|h|}{|x|}\,|x|^{-1/10},
		\]
		so Proposition~\ref{P-Holder} gives $(\mathrm H_2)$ whenever
		$\delta=3/5<\varepsilon$. Since
		$\max\{\frac23,\frac{10}{19}\}<\frac{10}{13}<1$ and $(\mathrm H)$ holds,
		(i) follows from Theorem~\ref{T-HpLq} and Examples~\ref{Ex-radialweights}
		and~\ref{Ex-symm}, while (ii) is the case $\varepsilon\le\delta-1/q=1/5$ of
		Remark~\ref{R-counterexample}.
	\end{example}
	
	\section*{Acknowledgement}
	We express our gratitude to Pablo Rocha (Universidad Nacional del Sur) and Marta S. Urciuolo (Universidad Nacional de C\'ordoba) for their helpful clarifications throughout the course of this work.
	
	\section*{AI Disclosure}
	During the preparation of this work, the authors used Claude ($Opus \ 5.5, \ Anthropic$) to assist in editing the manuscript and in writing parts of it, to check some technical steps, and to consult and check the bibliography. It was also used to construct some of the examples and auxiliary results, in particular those in Section~\ref{S-examples}.
	
	The main ideas and results of this work were developed by the authors, who verified all the arguments and proofs and take full responsibility for the content of the publication.

\end{document}